\documentclass{amsart}
\usepackage[dvipsnames]{xcolor}
\definecolor{darkgreen}{rgb}{0,0.5,0}
\usepackage{url}
\usepackage[
        colorlinks, citecolor=darkgreen,
        backref,
        pdfauthor={Isabel Rendell},
        pdftitle={Quadratic Chabauty for Atkin-Lehner quotients of modular curves via weakly holomorphic modular forms: Hodge Filtrations}
]{hyperref} 

\usepackage{tikz,enumerate,caption,amsfonts,amssymb,systeme,comment}
\usetikzlibrary{matrix,positioning,arrows,shapes,chains,calc,automata}
\usetikzlibrary{decorations.pathreplacing}
\usepackage[T1]{fontenc}
\usepackage[utf8]{inputenc}
\usepackage{cancel}
\usepackage{mathtools}
\usepackage{xcolor}
\usepackage{tikz-cd}
\usepackage{amssymb}
\usepackage{lipsum}

\newlength{\proofmargin}
\renewcommand*{\backref}[1]{}
\renewcommand*{\backrefalt}[4]{%
  \ifcase #1 %
    \relax
  \or
    $\uparrow$#2.%
  \else
    $\uparrow$#2.%
  \fi%
}

  \renewenvironment{thebibliography}[1]{%
    \begin{oldthebibliography}{#1}%
      \setlength{\parskip}{0ex}%
      \setlength{\itemsep}{0.1ex}%
      \setlength{\labelwidth}{1.4cm} 
  }%
  {%
    \end{oldthebibliography}%
  }

\DeclareFontFamily{U}{wncy}{}
\DeclareFontShape{U}{wncy}{m}{n}{<->wncyr10}{}
\DeclareSymbolFont{mcy}{U}{wncy}{m}{n}
\DeclareMathSymbol{\Sha}{\mathord}{mcy}{"58}

\DeclareMathOperator{\Q}{\mathbf{Q}}

\DeclareMathOperator{\Z}{\mathbf{Z}}
\DeclareMathOperator{\C}{\mathbf{C}}

\DeclareMathOperator{\dR}{\mbox{\scriptsize dR}}

\DeclareMathOperator{\End}{\mathrm{End}}

\newcommand{\Ker}{\mathrm{Ker}}

\newtheorem{theorem}{Theorem}[section]
\newtheorem{prop}[theorem]{Proposition}

\newtheorem{lemma}[theorem]{Lemma}
\newtheorem{corollary}[theorem]{Corollary}

\newtheorem{Definition}[theorem]{Definition}

\newtheorem{Remark}[theorem]{Remark}

\newtheorem{algorithm}[theorem]{Algorithm}

\makeatletter
\renewenvironment{proof}[1][\proofname]%
{%
\par\pushQED{\qed}\normalfont\topsep6\p@\@plus6\p@\relax%
\begin{list}{}{\rightmargin=8pt\leftmargin=\proofmargin}%
  \item[\hskip\labelsep\bfseries#1\@addpunct{.}]\ignorespaces
}{%
\popQED\end{list}\@endpefalse%
}%
\makeatother

\title[Quadratic Chabauty via weakly holomorphic modular forms: Hodge filtrations]{Quadratic Chabauty for Atkin-Lehner quotients of modular curves via weakly holomorphic modular forms: Hodge Filtrations}
\author{Isabel Rendell}
\address{Isabel Rendell, Department of Mathematics, King’s College London, Strand, London WC2R 2LS, UK}
\email{isabel.rendell@kcl.ac.uk}

\date{}

\begin{document}

\begin{abstract} For Atkin-Lehner quotients $X_0^+(N)$, of prime level and of genus at least 2, we provide an algorithm for computing one of the main objects in the quadratic Chabauty algorithm in terms of weakly holomorphic modular forms associated to the curve. In particular, the algorithm computes a Hodge filtration on a certain unipotent vector bundle with connection related to $X_0^+(N)$, which is crucial in computing the $p$-adic height which is used to define the finite set of $p$-adic points containing the rational points on $X_0^+(N)$. This improves the current Hodge filtration algorithm by replacing the input of an explicit plane model of the curve with weakly holomorphic modular forms to produce a faster computation. We implement our algorithm on the genus 7 modular curve $X_0^+(193)$, and discover congruences between iterated integrals of weight 2 cusp forms in the plus eigenspace for the Atkin-Lehner involution and single integrals of weight 2 cusp forms in the minus eigenspace for the Atkin-Lehner involution.\end{abstract}
\maketitle
\tableofcontents


\section{Introduction} \label{sec: intro}

\par The question of determining the rational points on the Atkin-Lehner quotient of the modular curve $X_0(N)$, $X_0^+(N) \coloneqq X_0(N) / w_N$, where $w_N$ is the Atkin-Lehner involution, is motivated by Mazur's Program B \cite{mazur2006rational}, which roughly amounts to determining the rational points on all modular curves. The non-cuspidal points on $X_0^+(N)$ correspond to unordered pairs of elliptic curves, along with a cyclic isogeny of degree $N$ between them, and the Atkin-Lehner involution, $w_N$, sends an isogeny to its dual. These non-cuspidal rational points can be divided into two types: those that correspond to elliptic curves with CM, or those that correspond to quadratic $\Q$-curves that are $N$-isogenous to their conjugates and points of this second type are known as \textit{exceptional}. The set of rational points is known for the curves $X_0^+(N)$ for $N$ prime such that the genus of $X_0^+(N)$ is between 2 and 6 due to \cite{balakrishnan2023quadratic} \cite{advzaga2021quadratic}. Furthermore, Galbraith conjectures \cite{galbraith2002rational} that for a prime $N >\!\!> 0$ the curve $X_0^+(N)$ has no exceptional rational points, and therefore for a prime $N >\!\!> 0$ there are no quadratic $\Q$-curves that are $N$-isogenous to their conjugates. In these cases, the rational points on $X_0^+(N)$ have been computed using the quadratic Chabauty algorithm, created by Balakrishnan-Dogra-M\"uller-Tuitman-Vonk \cite{balakrishnan2023quadratic}. The quadratic Chabauty algorithm computes the set of rational points for a nice curve of genus at least 2 satisfying some extra conditions, where a \textit{nice} curve is one that is smooth, projective, and geometrically irreducible. The theory behind quadratic Chabauty is a special case of Kim's generalisation \cite{kim2005motivic} of the Chabauty-Coleman method. 

\par For a nice curve $X$ of genus $g$, let $p$ be a prime of good reduction for $X$, let $J$ be its Jacobian, and let $r$ denote the rank of the rational points of $J$, so $r :=  \text{rank}J(\Q)$; this notation will remain throughout this paper. Before Faltings' proof of the Mordell Conjecture \cite{faltings1983endlichkeitssatze}, Chabauty \cite{chabauty1941points} asserted a special case: where $X$ satisfies $r < g$. The condition $r < g$ is known as the \textit{Chabauty-Coleman criterion}. Chabauty's proof involved embedding a curve $X$ into its Jacobian, $J$, and proving that the intersection $X(\Q_p) \cap \overline{J(\Q)}$ is finite, and thus $X(\Q)$ is finite. Coleman then adapted Chabauty's proof \cite{coleman1985effective} to give an explicit upper bound for the size of $X(\Q)$, still requiring that the Chabauty-Coleman criterion is satisfied. Coleman did this by constructing $p$-adic line integrals called \textit{Coleman integrals}, so that the rational points on $X$ lie in the kernel of these integrals. These integrals have finitely many zeroes, and then bounding the number of zeroes of these integrals gives an explicit upper bound on the number of rational points on $X$. Moreover, we can compute the kernel of these integrals, and thus obtain a finite set of $p$-adic points on $X$, containing $X(\Q)$, and then extract $X(\Q)$ from it. Therefore when the Chabauty-Coleman criterion is satisfied one can explicitly compute the rational points of a nice curve, by computing the vanishing locus of this Coleman integral; this process is known as the Chabauty-Coleman method. Unfortunately, the condition that the Chabauty-Coleman criterion is satisfied is very restrictive, and cannot be applied to compute rational points on curves with $r \geq g$.

\par In 2005, Kim \cite{kim2005motivic} produced a vast generalisation of the Chabauty-Coleman method by taking a unipotent quotient of the $\Q_p$-\'etale fundamental group of a curve $X$ and defining a corresponding subset of $p$-adic points. This subset of $p$-adic points is constructed so that it contains the rational points in order to obtain a \textit{nonabelian Chabauty} method. In particular, Kim's method produces an infinite tower of $p$-adic subsets, 
\begin{equation*}
    X(\Q) \subseteq ... \subseteq X(\Q_p)_n \subseteq ... \subseteq X(\Q_p)_2 \subseteq X(\Q_p)_1 \subseteq X(\Q_p),
\end{equation*}
so that the sets, $X(\Q_p)_n$, for all $n \in \mathbb{N}$, known as the \textit{Chabauty-Kim set} of \textit{depth} $n$, all contain $X(\Q)$. The depth one case, which concerns the set $X(\Q_p)_1$, roughly corresponds to the finite set of $p$-adic points produced in the Chabauty-Coleman method, provided the Chabauty-Coleman criterion is satisfied. The Chabauty-Kim method hopes to compute the set $X(\Q)$ by first explicitly computing the Chabauty-Kim sets $X(\Q_p)_n$, and then extracting the set $X(\Q)$ from the Chabauty-Kim sets. In particular, since $X(\Q_p)_n \subseteq X(\Q_p)_1$ for all $n \in \mathbb{N}$, the Chabauty-Kim method hopes to compute $X(\Q_p)_n$ for curves where $r \geq g$ by computing these smaller sets $X(\Q_p)_n$ for $n > 1$. For these such curves the Chabauty-Coleman criterion fails and we likely have that the set $X(\Q_p)_1$ is equal to the full set of $p$-adic points on $X$, $X(\Q_p)$, in which case $X(\Q_p)_1$ is no longer finite. For the process of computing $X(\Q_p)_n$ for some $n$ and extracting $X(\Q)$ to make sense we require that there is an integer $n$ such that the set $X(\Q_p)_n$ is finite and can be computed in practice. Kim conjectures that when $n >\!\!> 0$, $X(\Q_p)_n$ is finite \cite{kim2009unipotent}, which is implied by a special case of the Bloch-Kato conjectures, and further conjectures that when $n >\!\!> 0$, $X(\Q_p)_n = X(\Q)$ \cite{balakrishnan2018non}. In practice these sets are very difficult to compute because of the complexity of the nonabelian objects involved. However, successful instances of carrying out the Chabauty-Kim method relate to calculating the depth two set $X(\Q_p)_2$; this is the process known as quadratic Chabauty.

\par A critical development, due to the work of Balakrishnan-Dogra \cite{Balakrishnan_2018}, is the construction of a finite set of $p$-adic points containing $X(\Q_p)_2$ which is defined in terms of $p$-adic height functions, provided that $X$ satisfies certain conditions. The benefit of this is that these $p$-adic height functions can be computed explicitly. This extends previous work of Balakrishnan-Besser \cite{balakrishnan2015coleman} of using $p$-adic heights to determine integral points on an elliptic curve of rank 1 and, more generally, on odd degree hyperelliptic curves due to Balakrishnan-Besser-M\"uller \cite{balakrishnan2016quadratic}. Moreover, in the event that the set $X(\Q_p)_2$ is strictly greater than the set of rational points, one can employ techniques such as the Mordell-Weil sieve to show that given $p$-adic points do not come from a rational point. One of the first instances of carrying out quadratic Chabauty to compute the rational points on a modular curve is the split Cartan modular curve of level 13, $X_s^+(13)$, due to Balakrishnan-Dogra-M\"uller-Tuitman-Vonk \cite{balakrishnan2019explicit}. These authors then develop the quadratic Chabauty algorithm in \cite{balakrishnan2023quadratic}, for a modular curve $X$ with $g \geq 2$ such that $r = g$, and $\rho(J) > 1$, where $\rho(J)$ is the N\'eron-Severi rank of the Jacobian. This requires an explicit plane model of the curve as an input and thus we will sometimes refer to this as the \textit{model} algorithm. In addition to being implemented on Atkin-Lehner quotients $X_0^+(N)$ for $N$ prime such that the genus is between 2 and 6 \cite{balakrishnan2023quadratic}, \cite{balakrishnan2021two}, \cite{advzaga2021quadratic}, the model algorithm has been used to compute the rational points on the modular curves $X_{S_4}(13), X_{ns}^+(17)$ \cite{balakrishnan2023quadratic}, and extended to number fields \cite{balakrishnan2025rationalpointsnonsplitcartan}. For all of these examples, the Chabauty-Coleman criterion fails since $r = g$.

\par Therefore, the two aims of the quadratic Chabauty method for a curve $X$ are first to show that the set $X(\Q_p)_2$ is finite, and second to construct explicit functions on $X(\Q_p)$ which cut out (a finite set containing) $X(\Q_p)_2$, and hence containing $X(\Q)$. For the first step, Balakrishnan-Dogra assert that when $r < \rho(J) + g - 1$, the set $X(\Q_p)_2$ is finite. Since we have $\rho(J) \geq 1$, this condition is always at least as strong as the Chabauty-Coleman criterion. It turns out that modular curves are a natural candidate for quadratic Chabauty, as they often satisfy this condition, as shown by Siksek \cite{siksek2017quadratic}, and so the set $X(\Q_p)_2$ is often finite for modular curves. Additionally, Dogra-Le Fourn \cite{dogra2021quadratic} extend this to quotients of modular curves to assert that the set $X(\Q_p)_2$ is finite for $X = X_0^+(N), X_{ns}^+(N)$ for $N$ prime such that the genus is at least 2. 

\par We now discuss the second aim of quadratic Chabauty: to construct explicit functions on $X(\Q_p)$ which cut out (a finite set containing) $X(\Q_p)_2$. Balakrishnan-Dogra have proven that these functions can be written in terms of $p$-adic heights. A fundamental property of the global $p$-adic height function, $h$, is that it has a decomposition into a sum of local factors $h_v$ for all primes $v$ of $\Q$. In particular, we have $h = h_p + \sum_{v \neq p} h_v$, where $h_p$ is known as the local $p$-adic height at $p$, and $h_v, v \neq p$ are called the local $p$-adic heights away from $p$. It will be explained in §\ref{sec: prelims}, that the local $p$-adic height at $p$, $h_p$, is defined via a function that associates to each $\Q_p$-point on $X$, $x$, a filtered $\phi$-module (or mixed extension of filtered $\phi$-modules), denoted by $D_{cris}(A_Z(b, x))$, where $D_{cris}$ is Fontaine's $D_{cris}$-functor, and $A_Z(b, x)$ is a certain Galois representation associated to a point $x$ on $X$. Once $D_{cris}(A_Z(b, x))$ is computed, it is straightforward to compute the local $p$-adic height function. In particular, to compute $D_{cris}(A_Z(b, x))$, it is sufficient to compute a Hodge filtration and Frobenius structure on a certain unipotent vector bundle with connection on $X$, and this unipotent vector bundle is uniquely determined by the entries of a $(2g + 2) \times(2g+2)$ matrix with entries in $\text{H}^1_{\dR}(X/\Q)$. In turn, computing the Hodge filtration and Frobenius structure is equivalent to computing two matrices $\lambda^{\text{Fil}}, \lambda^\phi$, respectively, and computing the matrix $\lambda^{\text{Fil}}$ in terms of weakly holomorphic modular forms is exactly the subject of this paper.

\par The drawback of the current quadratic Chabauty algorithm is that it requires an explicit plane model of the curve as an input. These models are inefficient to compute with since the coefficients of such a model grow rapidly with the genus of the modular curve. All of the past examples of curves that the model algorithm has been implemented on have a genus of at most 6. This provides motivation for the creation of a more efficient algorithm, which removes the need for computing with the model of the curve. To do this, we would like to exploit the fact that we have a modular curve and replace this input of a model with modular forms associated with the curve. We refer to such an algorithms and methods as \textit{model-free}. Another disadvantage of the model-dependent approach is that it obscures interpretations of quadratic Chabauty in terms of modular forms or the arithmetic of modular curves. By comparison, our method, which uses weakly holomorphic modular forms, uncovers congruences between single and iterated integrals on $X_0^+(N)$, as we see in \S \ref{sec: example}. Existing model-free methods include computing single Coleman integrals on modular curves due to Chen-Kedlaya-Lau \cite{chen2024colemanintegrationmodularcurves}, and computing the global $p$-adic height function in terms of special values of $p$-adic L-functions due to Hashimoto \cite{Hashimoto_2023}.

\par From now on, let $X$ be an Atkin-Lehner quotient $X_0^+(N)$ for $N$ prime. Since the defining equations for the (finite set containing) $X(\Q_p)_2$ can be defined in terms of $p$-adic heights, in order to develop a model-free algorithm, one could compute $p$-adic height functions in terms of modular forms associated to $X$. Indeed, the main result of this paper is an algorithm to compute a Hodge filtration in terms of modular forms associated to $X$, Algorithm \ref{alg: model-free HF}, which is the first step towards a model-free quadratic Chabauty algorithm on Atkin-Lehner quotients of prime level. An initial step in the quadratic Chabauty algorithm is to compute a basis of $\text{H}_{\dR}^1(X/\Q)$, which is usually defined in terms of the model of the curve. Therefore, one of the contributions of this paper is a method to compute a basis of $\text{H}_{\dR}^1(X/\Q)$ in terms of $q$-expansions of cusp forms, which will be the focus of \S \ref{subsec: computing H1dR}. Once we have such a basis, we explain in \S \ref{subsec: calculating hodge} how we then obtain the matrix $\lambda^{\text{Fil}}$ in terms of weakly holomorphic modular forms, which includes the use of the Hecke operator $T_p$ for a certain prime $p$. Therefore, we produce the following algorithm, which we have implemented in \texttt{Magma} \cite{MR1484478}, and which is later given as Algorithm \ref{alg: model-free HF}.

\begin{algorithm}[Model-free Computation of the Hodge filtration for $X_0^+(N)$ for $N$ prime] \label{alg: model free HF in intro}
    \hfill \break 
    Input:
    \begin{itemize}
        \item A prime $N$ such that $X_0^+(N)/\Q$ is a modular curve such that $g > 1$.
        \item A prime $p$ such that the Hecke operator $T_p$ does not lie in $\Z \subset \text{End}(J)$.
    \end{itemize}
    Output: the matrix $\lambda^{\text{Fil}}$, defined in terms of weakly holomorphic modular forms.
\end{algorithm}

\par We will see that we obtain one half of the basis of $\text{H}_{\dR}^1(X/\Q)$ in terms of $q$-expansions of weight 2 cusp forms for $X$, which in turn defines a basis of $q$-expansions of $\text{H}^0(X, \Omega^1)$, and we denote such a basis by $\{ \omega_0, ..., \omega_{g-1}\}$. For the second half of the basis we compute a function $f_{\dR}$ in terms of weight 12 cusp forms so that $\{\omega_0, ..., \omega_{g-1}, f_{\dR} \cdot \omega_0, ..., f_{\dR} \cdot \omega_{g-1}\}$ is a basis of $q$-expansions of $\text{H}_{\dR}^1(X/\Q)$, and we can verify that it is a basis using Serre's formula. Moreover, we describe in detail how to implement our algorithm on the modular curves $X_0^+(67)$ and $X_0^+(193)$, which have genus 2 and 7 respectively. We also discover congruences between iterated integrals of weight 2 cusp forms in the plus eigenspace for the Atkin-Lehner involution and single integrals of weight 2 cusp forms in the minus eigenspace for the Atkin-Lehner involution by computing matrices of coefficients of the $q$-expansions of these objects. 



\subsection{Outline}
\par The structure of the paper is as follows. In §\ref{sec: prelims} we introduce the necessary background for the current quadratic Chabauty algorithm, by giving a brief outline of the Chabauty-Kim method in \S \ref{subsec: C-K method}, introducing $p$-adic height functions in \S \ref{subsec: QC}, and then outlining the current algorithm for computing the Hodge filtration in \S \ref{subsec: hodge}. In §\ref{sec: model-free theory} we first explain how to construct a basis of $\text{H}_{\dR}^1(X/\Q)$ in terms of modular forms related to our modular curve in \S \ref{subsec: computing H1dR}, and then explain how to compute the Hodge filtration in terms of modular forms in \S \ref{subsec: calculating hodge}. Finally, in \S \ref{sec: example} we implement our algorithm on a genus 2 curve $X_0^+(67)$ and a genus 7 curve $X_0^+(193)$, and investigate congruences between integrals by computing the matrices of coefficients of cusp forms for $X_0^+(67)$ and the genus 3 curve $X_0^+(97)$.

\subsection*{Acknowledgments}  I am extremely grateful to my supervisor, Netan Dogra, for suggesting this project, for very many valuable conversations, and for insightful comments on earlier drafts. This work was supported by the Engineering and Physical Sciences Research Council [EP/S021590/1]. The EPSRC Centre for Doctoral Training in Geometry and Number Theory (The London School of Geometry and Number Theory), University College London. This research was conducted at King's College London.

\section{Quadratic Chabauty: theory and algorithms}
\label{sec: prelims}



\par In this section, we provide an overview of the background on the theory of quadratic Chabauty, $p$-adic heights and the model algorithm of computing the Hodge filtration. Since the subject of this paper is a model-free algorithm for computing the Hodge filtration, the focus of this section will be explicitly computing the objects involved in quadratic Chabauty, especially focusing on the existing algorithm to compute the Hodge filtration in \S \ref{subsec: hodge}, as this will be required in \S \ref{sec: model-free theory}. Therefore, we only include a brief overview on the theory underpinning quadratic Chabauty, the Chabauty-Kim method, at the start of this section.

\subsection{The Chabauty-Kim Method} \label{subsec: C-K method}

\par Since quadratic Chabauty is a successful implementation of the Chabauty-Kim method, we begin \S \ref{sec: prelims} by providing a very brief overview of Kim's generalisation of the Chabauty-Coleman method. This will be a deliberately brief overview, the sole purpose of which is to provide some insight into how the subset of $p$-adic points, $X(\Q_p)_2$, is constructed and the conditions which guarantee that this set $X(\Q_p)_2$ is finite. In \S \ref{sec: intro} we saw how the Chabauty-Coleman method relies on associating the Jacobian of a curve in order to determine the rational points. This is achieved by constructing a subset of $p$-adic points, containing the rationals, providing conditions such that this set is finite, and then explicitly computing this finite subset. For Kim's generalisation, rather than associate a curve's Jacobian, he employs a unipotent quotient of the $\Q_p$-\'etale fundamental group of the curve, which is motivated by the section conjecture as follows. 

\par In \cite{grothendieck1997brief}, Grothendieck relates rational points on a curve $X$ to the geometric \'etale fundamental group $\pi_1^{\text{\' et}}(\overline{X}, b)$ of $X$ with basepoint $b$. Specifically, he conjectures that the map 
\begin{align*}
    X(\mathbb{Q})  &\longrightarrow \text{H}^1(G_{\mathbb{Q}},  \pi_1^{\text{\'et}}(\overline{X}, b)), \\
    x &\mapsto [\pi_1^{\text{\'et}}(\overline{X}; b, x)],
\end{align*}
which attaches to each rational point the class of the Galois representation defined by the corresponding path torsor of the algebraic fundamental group, should be an isomorphism. Therefore, this conjectured isomorphism provides a potential new avenue through which we can study rational points on a curve $X$. However, it seems that taking the full set of such torsors has too little structure in order to be utilised for explicit computations. Therefore, we are motivated to take a quotient of $\pi_1^{\text{\' et}}(\overline{X}, b)$ that has a suitable amount of structure, so that we can work with more classes of curves than the Chabauty-Coleman method. Hence, the Chabauty-Kim method relies on choosing such a quotient and defining a corresponding set of $p$-adic points containing the rationals. The hope is then that this set can be proven to be finite in certain circumstances and that it can also be explicitly computed. We now discuss the construction of this set of $p$-adic points and conditions which guarantee it is finite. Then, the remainder of \S \ref{sec: prelims} will focus on how, under certain circumstances, it can indeed by explicitly computed. For the remainder of this paper, we are interested in the case that our curve does not satisfy the Chabauty-Coleman criterion, so that we have $r \geq g$.

\par Assume that $X(\Q) \neq \emptyset$ and fix $b \in X(\Q)$. Let $p$ be a prime of good reduction, let $T_0$ denote the set of bad primes for $X$, and we set $T \coloneqq T_0 \cup \{p\}$. Let $G_T$ be the maximal quotient of $G_{\Q}$ unramified outside $T$. Let $U$ be a finite dimensional Galois-stable unipotent quotient of $\pi^{\text{\'et}}_1(X_{\overline{\Q}, b})_{\Q_p}$, the $\Q_p$-pro-unipotent completion of \'etale fundamental group of $X_{\overline{\Q}}$ with base point $b$. The following commutative diagram is vital to the Chabauty-Kim method
\begin{equation} \label{eqn: C-K commutative diagram}
    \begin{tikzcd}
    X(\Q) \arrow{d}{j_U}  \arrow[swap]{r} & X(\Q_p) \arrow{d}{j_{U, p}}  \\ 
     \text{H}^1(G_T, U) \arrow{r}{\text{loc}_{U, p}} &  \text{H}^1(G_p, U), 
  \end{tikzcd}
\end{equation}
where $j_U, j_{U, p}$ denote the global and local unipotent Kummer maps defined in \cite{kim2005motivic}, \cite{kim2009unipotent}, respectively. Moreover, for $v \neq p$, we also have local analogous maps $j_{U, v}: X(\Q_v) \mapsto \text{H}^1(G_v, U)$. A crucial result of Kim \cite{kim2005motivic}, \cite{kim2009unipotent} is that the nonabelian pointed continuous cohomology sets $\text{H}^1(G_T, U), \text{H}^1(G_v, U)$ for all $v$ are affine algebraic varieties over $\Q_p$. Kim also shows that the localisation maps, $\text{loc}_{U, v}$, for all $v$, are algebraic and that the crystalline torsors have the structure of the $\Q_p$-points on a subvariety. Moreover, we can obtain another commutative diagram by replacing the bottom row of (\ref{eqn: C-K commutative diagram}) with a certain subscheme and subvariety, respectively, which contain the images of the rational points under the maps $j_U, j_{U, p}$, and this is done in terms of crystalline classes. The benefit of dealing with the certain subscheme and subvariety is the hope that these are more amenable to explicit computation. We now define the certain subscheme of $\text{H}^1(G_T, U)$.

\begin{Definition}[\cite{Balakrishnan_2018} Definition 2.2]
    We define the Selmer variety $\text{Sel}(U)$ to be the reduced scheme associated to the subscheme of $\text{H}^1(G_T, U)$ containing those classes $c$ such that 
    \begin{itemize}
        \item $\text{loc}_{U, p}(c)$ is crystalline
        \item $\text{loc}_{U, v}(c) \in j_{U, v}(X(\Q_v))$ for all $v \neq p$
        \item the projection of $c$ to $\text{H}^1(G_T, U)$ comes from an element of $J(\Q) \otimes \Q_p$.
    \end{itemize}
\end{Definition}

\par For the subvariety of  $\text{H}^1(G_p, U)$, Kim \cite{kim2009unipotent} proves that the local Bloch-Kato Selmer scheme of crystalline classes, denoted by $\text{H}^1_f(G_p, U)$, is indeed a subvariety of $\text{H}^1(G_p, U)$. Moreover, the fact that $\text{H}^1_f(G_p, U)$ contains the image of $X(\Q_p)$ under $j_{U, p}$ is due to a result of Olsson (\cite{olsson2011towards} Theorem 1.4), so that we have $j_{U, p}(X(\Q_p)) \subset \text{H}^1_f(G_p, U)$. Therefore, we obtain another commutative diagram by replacing $\text{H}^1(G_T, U), \text{H}^1(G_p, U)$ in diagram (\ref{eqn: C-K commutative diagram}) with $Sel(U), \text{H}^1_f(G_p, U)$ respectively. As a result, we can construct a set of $p$-adic points containing the rationals in terms of $U$ in order to obtain the following tower of sets
\begin{equation*}
    X(\Q) \subset X(\Q_p)_U \coloneqq j_{U, p}^{-1}(\text{loc}_{U, p}(\text{Sel}(U))) \subset X(\Q_p).
\end{equation*}
Now suppose that $U$ is a Galois-stable quotient of the maximal $n$-unipotent quotient of $\pi^{\text{\'et}}_1(X_{\overline{\Q}, b})_{\Q_p}$, which we denote by $U_n$. Then we obtain the following definition of the Chabauty-Kim sets, $X(\Q_p)_n$,
\begin{equation} \label{X(Q_p)_2 subset X(Qp)_U}
    X(\Q) \subset X(\Q_p)_n \coloneqq X(\Q_p)_{U_n} \subset X(\Q_p)_U.
\end{equation}

\par Therefore, we have constructed a set of $p$-adic points containing the set $X(\Q)$. In order for computing the set $X(\Q_p)_n$ to make sense as a method to compute the rational points, we must first have that $X(\Q_p)_n$ is finite and secondly that it is computable in practice. As mentioned, for the first of these conditions, Kim's conjectures that we have $n >\!\!> 0$, $X(\Q_p)_n$ is finite \cite{kim2009unipotent}. Moreover, he later conjectures that we eventually obtain the rational points themselves, so that for $n >\!\!> 0$, $X(\Q_p)_n = X(\Q)$ \cite{balakrishnan2018non}. In addition to these conjectures, Kim proves the following result regarding finiteness of $X(\Q_p)_U$.

\begin{theorem}[Kim \cite{kim2009unipotent}] \label{thm: X(Qp)_U is finite}
    Suppose that $\text{loc}_{U, p}$ is non-dominant. Then $X(\Q_p)_U$ is finite.
\end{theorem}
\par For the remainder of \S \ref{subsec: C-K method} we will give finiteness results for the Chabauty-Kim set $X(\Q_p)_2$ which is foundational to the quadratic Chabauty method as we will see now. Suppose that the set $X(\Q_p)_1$ cut out by classical Chabauty-Coleman is infinite. The goal of the quadratic Chabauty method is to:
\begin{enumerate}[(i)] 
    \item show that $X(\Q_p)_2$ is finite,
    \item construct explicit functions on $X(\Q_p)$ cutting out a finite set containing $X(\Q_p)_2$.
\end{enumerate}
The remainder of this section will tackle (i), leaving (ii) for §\ref{subsec: QC}. To state the fundamental result regarding (i), we first need to recall the following definition. 

\begin{Definition}[N\'eron-Severi group] \label{def: rho(J)}
For a variety $X$, define the N\'eron-Severi group of $X$, $\text{NS}(X)$, as the group of divisors modulo algebraic equivalence. Therefore, we have $\text{NS}(X) = Pic(X)/Pic^0(X)$. The N\'eron-Severi rank of $X$, $\rho(X)$, is the rank of $\text{NS}(X)$, and is often referred to as the Picard number.
\end{Definition}

\par We can now give the fundamental result concerning finiteness of $X(\Q_p)_2$, which was proven via Theorem \ref{thm: X(Qp)_U is finite} and by computing and bounding dimensions of $\text{H}_f^1(G_T, U)$, and $\text{H}_f^1(G_p, U)$.

\begin{theorem}[\cite{Balakrishnan_2018}] \label{thm: when is X(Qp)2 finite}
    Suppose that 
    \begin{equation} \label{QC condition}
        r < g + \rho(J) - 1,
    \end{equation}
    then $X(\Q_p)_2$ is finite.
\end{theorem}

\begin{Remark}
    We have $\rho(J) \geq 1$, and so this condition is always at least as strong as the Chabauty-Coleman condition that $r < g$, and therefore can be applied in more cases.
\end{Remark}

\par We call (\ref{QC condition}) the \textit{quadratic Chabauty condition}. As we have already mentioned, Siksek showed that the quadratic Chabauty condition is satisfied frequently by modular curves in \cite{siksek2017quadratic}. Furthermore,  work of Dogra and Le Fourn \cite{dogra2021quadratic} gives finiteness of $X(\Q_p)_2$ for quotients of modular curves in many circumstances.

\begin{theorem}[\cite{dogra2021quadratic}] \label{thm: finiteness for X0+(N)}
\hfill  

\begin{itemize}
        \item For all prime $N$ such that $g(X_0^+(N)) \geq 2$, $X_0^+(N)(\Q_p)_2$ is finite for any $p \neq N$.
        \item For all prime $N$ such that $g(X_{ns}^+(N)) \geq 2$ and $X_{ns}^+(N)(\Q) \neq \emptyset$, $X_{ns}^+(N)(\Q_p)_2$ is finite for any $p \neq N$.
    \end{itemize}
\end{theorem}

\par Therefore, we have finiteness of $X(\Q_p)_2$ in many cases involving modular curves and their quotients, and so it makes sense to try to apply quadratic Chabauty to these curves to obtain $X(\Q)$. In particular, Theorem \ref{thm: finiteness for X0+(N)} gives that we always have finiteness of $X(\Q_p)_2$ for the subject of this paper: $X_0^+(N)$ for $N \neq p$ prime and $g > 1$.

\subsection{Quadratic Chabauty for rational points: $p$-adic heights} \label{subsec: QC}

 \par In general, the objects involved in the Chabauty-Kim method are very hard to explicitly compute. One approach which has been highly successful in implementation is the method of quadratic Chabauty. The success of this is due to relating the functions involved in Kim's method to $p$-adic heights. Suppose $r = g$ (and so we don't satisfy the conditions for Chabauty-Coleman), and that the $p$-adic closure of $J(\Q)$ has finite index in $J(\Q_p)$. Let $\text{AJ}_b$ be the Abel-Jacobi morphism, which is the map
 \begin{align*}
     \text{AJ}_b: X(\Q_p) &\rightarrow \text{H}^0(X_{\Q_p}, \Omega^1)^* \\
     x &\mapsto \left(\omega \mapsto \int_b^x \omega \right).
 \end{align*}
Under our assumptions, $\text{AJ}_b$ induces an isomorphism $J(\Q) \otimes \Q_p \cong \text{H}^0(X_{\Q_p}, \Omega^1)^*$, and we call this isomorphism $\alpha : J(\Q) \otimes \Q_p \rightarrow \text{H}^0(X_{\Q_p}, \Omega^1)^*$. Since we have an isomorphism, $\alpha$, global points cannot be cut out of local points using linear relations in the Abel-Jacobi map, which highlights the significance of the condition that $r = g$. Therefore, the idea of the quadratic Chabauty method is to replace linear relations with \textit{bilinear} relations. Suppose we can find a function $\theta : X(\Q_p) \rightarrow \Q_p$ 
 and a finite set $\Upsilon \subset \Q_p$ with the following properties:
 \begin{itemize}
     \item On each residue disc $]x[ \subset X(\Q_p)$, the map
     \begin{equation*}
         (\text{AJ}_b, \theta) : X(\Q_p) \longrightarrow \text{H}^0(X_{\Q_p}, \Omega^1)^* \times \Q_p 
     \end{equation*}
     has Zariski dense image and is given by a convergent power series. 
     \item There exists an endomorphism $E$ of $\text{H}^0(X_{\Q_p}, \Omega^1)^*$, a functional $c \in \text{H}^0(X_{\Q_p}, \Omega^1)^*$, and a bilinear form $B : \text{H}^0(X_{\Q_p}, \Omega^1)^* \otimes \text{H}^0(X_{\Q_p}, \Omega^1)^* \rightarrow \Q_p$, such that, for all $x \in X(\Q)$, 
     \begin{equation} \label{theta eqn}
         \theta(x) - B(\text{AJ}_b(x), E(\text{AJ}_b(x)) + c) \in \Upsilon.
     \end{equation}
 \end{itemize}
 
\par These two properties produce a finite set of $p$-adic points containing $X(\Q)$ since the first property implies that only finitely many $p$-adic points can satisfy equation (\ref{theta eqn}), and all rational points satisfy it by the second property. Finiteness of this set is proven by local analytic information and global arithmetic information. We call $(\theta, \Upsilon)$ a \textit{quadratic Chabauty pair}, and the objects $E, c,$ and $B$ of a quadratic Chabauty pair will be referred to as its endomorphism, constant, and pairing, respectively. We now discuss how to compute the functions in (\ref{theta eqn}) in terms of $p$-adic heights.

\par The first implementation of quadratic Chabauty was computing integral points on rank 1 elliptic curves due to Balakrishnan-Besser \cite{balakrishnan2015coleman}, and on odd degree hyperelliptic curves due to Balakrishnan-Besser-M\"uller \cite{balakrishnan2016quadratic}. This is done via the global Coleman-Gross $p$-adic height function, $h$, and its decomposition into local $p$-adic heights, $h_v$, and utilising the different properties of the local $p$-adic height at $p$, $h_p$, and the local heights away from $p$, $h_v$, with $v \neq p$. When $x$ is a rational but not integral point, the values of $h_v(x)$ for $v \neq p$ cannot be controlled and so a different $p$-adic height is adopted in order to compute the set of rational points.
 
\par Balakrishnan-Dogra \cite{Balakrishnan_2018} show that in order to apply quadratic Chabauty to compute rational points, quadratic Chabauty pairs are defined by associating certain Galois representations to points of $X$, so that we can then use Nekov\'a\v r's $p$-adic height functions \cite{nekovar1990p} on these certain Galois representations. We first discuss how these Galois representations are constructed, and then how Nekov\'a\v r's $p$-adic height functions can be used to produce quadratic Chabauty pairs. The Galois representations are defined in terms of a nice correspondence $Z$ on $X \times X$, which is constructed by nontrivial elements of $\Ker(NS(J) \rightarrow NS(X))$. This kernel is a free abelian group of rank equal to $\rho(J) - 1$, and so such a correspondence always exists when $\rho(J) > 1$. Let $K \in \{\Q, \Q_p\}$ and let $G_K$ denote the absolute Galois group of $K$. We can attach to any such choice of nice correspondence, $Z$, a suitable quotient $U_Z$ of the $\Q_p$-pro-unipotent fundamental group of $X_{\overline{\Q}}$, which in turn defines a certain family of Galois representations. We denote such representations by $A_Z(b, x)$, for each $x$, so that we have a map
\begin{align*}
    X(K) &\longrightarrow \{G_K \rightarrow \text{GL}_{2g+2}(\Q_p) \} / \sim \\
    x &\mapsto A_Z(b, x).
\end{align*}
A more detailed explanation can be found in \cite{Balakrishnan_2018}.  Let $V = \text{H}_{\text{\'et}}^1(X_{\overline{K}}, \Q_p)^*$, and let $V_{\dR}$ be the image of $V$ under Fontaine's $D_{\text{cris}}$-functor, so that $V_{\dR} = D_{\text{cris}}(V) = \text{H}^1_{\dR}(X/\Q_p)^*$. In addition, let $\chi_p : G_K \rightarrow \Q_p^\times$ be the $p$-adic cyclotomic character, and let $\Q_p(1)$ be the one-dimensional representation given by $\chi_p$. With respect to a suitable choice of basis, the representation $A_Z(b, x)$ is lower triangular, of the form
\begin{equation*}
        g \in G_K \mapsto \begin{pmatrix}
            1 & 0 & 0 \\
            * & \rho_V(g) & 0 \\
            * & * & \chi_p(g)
        \end{pmatrix},
\end{equation*}
where $\rho_V$ corresponds to $V$. Representations of this form, which admit a $G_K$-stable filtration with graded pieces $\Q_p(1), V, \Q_p$, along with certain isomorphisms, are referred to as \textit{mixed extensions}, see \cite{balakrishnan2019explicit} Definition 3.4. Therefore, we can define a mixed extension for each point $x \in X(\Q)$, and Nekov\'a\v r's local $p$-adic heights are functions on mixed extensions as we will see now.

\par We now discuss Nekov\'a\v r's $p$-adic height functions. Denote Nekov\'a\v r's global $p$-adic height by $h$. Firstly, $h$ has a decomposition into a sum of local height functions, $h_v$, defined locally at every finite place $v$, so that we have:
\begin{equation*} \label{eqn: nek decomposition}
        h(A_Z(b, x)) = h_p(A_Z(b, x)) + \sum_{v \neq p} h_v(A_Z(b, x)).
    \end{equation*}
    These local height functions have the following important properties: 
    \begin{enumerate}[(A)]
        \item for $v = p$, the map $x \mapsto h_p(A_Z(b, x))$ extends to a locally analytic function $\theta : X(\Q_p) \rightarrow \Q_p$ by Nekov\'a\v r's construction \cite{nekovar1990p}.
        \item for $v \neq p$ the local heights $h_v(A_Z(b, x))$ take on a finite set of values $\Upsilon \subset \Q_p$ for $x \in X(\Q_v)$ by a result of Kim and Tamagawa \cite{kim2008component}. In particular, the value of $h_v(A_Z(b, x))$ is identically zero when $v$ is a prime of potential good reduction.
    \end{enumerate}
    Secondly, since local $p$-adic height functions are functions on mixed extensions, $h$ can also be viewed as a function on mixed extensions. Moreover, \cite{balakrishnan2019explicit} \S 3.1, explains that, roughly, these mixed extensions can be mapped to $\text{H}^1_f(G_{\Q}, V)$, and the global height $h$ factors through this map but the local $p$-adic height functions do not. As a result, $h$ can be viewed as a function on $\text{H}^1_f(G_{\Q}, V)$, and in fact Nekov\'a\v r proves in \cite{nekovar1990p} Theorem 4.11 that $h$ defines a bilinear pairing 
    \begin{equation*}
        h : \text{H}^1_f(G_{\Q}, V) \times \text{H}^1_f(G_{\Q}, V) \longrightarrow \Q_p. 
    \end{equation*}
    Furthermore, under our assumptions at the start of \S \ref{subsec: QC}, we have that $h$ defines a bilinear pairing on $\text{H}^0(X_{\Q_p}, \Omega^1)^*$ via our isomorphism $\alpha$. We will now explain how combining the fact that $h$ defines a bilinear function on $\text{H}^0(X_{\Q_p}, \Omega^1)^*$ and the properties of the local $p$-adic height functions means that $(\theta, \Upsilon)$ produces a quadratic Chabauty pair. In particular, we can build a function as in $(\ref{theta eqn})$ with $B = h \circ \alpha^{-1}$, and the endomorphism, $E$, is the one induced by $Z$, by \cite{balakrishnan2019explicit} \S 3.4. 

\par By (B) above and the decomposition of $h$, we know that there exists a finite set $\Upsilon \subset \Q_p$ such that for any $x \in X(\Q)$ we have  
\begin{equation} \label{finiteness: h - h_p in Up}
    h(A_Z(b, x)) - h_p(A_Z(b, x)) \in \Upsilon.
\end{equation}
Moreover, we have that $h$ defines a bilinear pairing on $\text{H}^0(X_{\Q_p}, \Omega^1)^*$ and by \cite{balakrishnan2023quadratic} §2.3, $h$ may be extended to a bilinear map on $X(\Q_p)$ so that we have $h : X(\Q_p) \rightarrow \Q_p$. In addition, by (A) above we know $x \mapsto h_p(A_Z(b, x))$ also extends to give a locally analytic function $\theta : X(\Q_p) \rightarrow \Q_p$. Therefore, we can define a map $\rho_Z \coloneqq h - h_p$, which we call the \textit{quadratic Chabauty function} so that we have a function
\begin{equation*}
    \rho_Z : X(\Q_p) \longrightarrow \Q_p,
\end{equation*}
which is known to be Zariski dense on every residue disc. Therefore, the fact that all rational points on $X$ are mapped into the finite set $\Upsilon$ under this map, as in (\ref{finiteness: h - h_p in Up}), recovers the fact that the set $X(\Q)$ is finite. Since the map $\rho_Z$ is defined on $X(\Q_p)$, it is defined on $X(\Q_p)_2$, and Balakrishnan-Dogra \cite{Balakrishnan_2018} show that every element of the Chabauty-Kim set $X(\Q_p)_2$ is also mapped into $\Upsilon$ under the map $\rho_Z$. Therefore, we obtain the following tower of finite sets, 
\begin{equation} \label{def: X(Qp)QC set}
    X(\Q) \subset X(\Q_p)_2 \subset \{ x \in X(\Q_p) : h(A_Z(b, x)) - h_p(A_Z(b, x)) \in \Upsilon \}.
\end{equation}
The importance of constructing this largest set is that it is defined in terms of $p$-adic height functions which are computable objects. Therefore, to compute the rational points of $X$, we can compute a $p$-adic approximation of the largest finite set, which by construction contains the rationals, and employ techniques such as the Mordell-Weil sieve in the event that this set is strictly larger than $X(\Q)$.

\par Now, suppose that we have enough rational points, $P_1, ..., .P_m$, that satisfy the conditions laid out in \S1.4 of \cite{balakrishnan2019explicit}. Then, by the discussion above, if we can compute the following $p$-adic heights below, then we have a method to explicitly compute a finite subset of $X(\Q_p)$, containing the set $X(\Q)$, as defined in (\ref{def: X(Qp)QC set}):
\begin{enumerate}[(a)]
    \item Determine the set of values that $h_v(A_Z(b, x))$ can take for $x \in X(\Q_v)$ and $v \neq p$. 
    \item Expand the function $x \mapsto h_p(A_Z(b, x))$ into a $p$-adic power series on every residue disc.
    \item Evaluate $h(A_Z(b, P_i))$ for $i = 1, ..., m$.
\end{enumerate}
\par In the event that our curve $X$ has potentially good reduction everywhere (for example $X_s^+(13)$), then all local heights away from $p$ are trivial by the properties of $h_v$ for $v \neq p$, (B). In this case the three steps above reduce to solving (b) due to the decomposition of $h$. In this paper, our algorithm is implemented on the modular curve $X_0^+(N)$, for $N$ prime, which has good reduction away from $N$. Therefore, when $v \neq p, N$, for all $x \in X(\Q_v)$, we have that $h_v(A_Z(b, x))$ is equal to 0. Despite the fact that $X_0^+(N)$ does not have potentially good reduction at $N$, we can use the following result from \cite{balakrishnan2023quadratic} which implies that when applying quadratic Chabauty to $X_0^+(N)$, $N$ prime, there are no nontrivial contributions to the height away from $p$. 

\begin{lemma}[\cite{balakrishnan2023quadratic}, Lemma 5.2] \label{lemma: all local heights away from p are trivial}
The local height $h_N$ is trivial on $X_0^+(N)(\Q_N)$.
\end{lemma}

\par As a result, there are no nontrivial local heights $h_v$ for all $v \neq p$, so that the 3 steps above reduce to computing the local height $h_p(A_Z(b, x))$ only. Therefore, creating a model-free quadratic Chabauty algorithm for the modular curves $X_0^+(N)$ for $N$ prime reduces to the problem of computing the local $p$-adic height, $h_p$, without first computing an explicit model of the curve. Now that we have established the central role the local $p$-adic height at $p$, $h_p$, plays in implementing quadratic Chabauty on the modular curves $X_0^+(N)$, we now focus on how to actually compute it. Let $D_{\text{cris}}(A_Z(b, x))$ denote the image of $A_Z(b, x)$ under Fontaine's $D_{\text{\text{cris}}}$-functor. Then $D_{\text{cris}}(A_Z(b, x))$ is a filtered $\phi$-module, which is a finite dimensional $\Q_p$ vector space $W$ equipped with an exhaustive and separated decreasing filtration $\text{Fil}^i$ and an automorphism $\phi = \phi_W$, as in \cite{balakrishnan2019explicit} Definition 3.3. By Nekov\'a\v r's definition, to construct $h_p(A_Z(b, x))$, it suffices to explicitly describe its Hodge filtration and its Frobenius structure.

\par A key result is that $D_{\text{cris}}(A_Z(b, x))$ can be described as the pullback along $x$ of a certain universal connection. We denote this universal connection by $\mathcal{A}_Z$, and we will give more details of how it is defined in \S \ref{subsec: hodge}. In particular, we have the following isomorphism of filtered $\phi$-modules from \cite{balakrishnan2019explicit} §5, stated as Lemma 5.16 in \cite{balakrishnan2023computation}. 
\begin{lemma} \label{lemma: eqn: D cris = x^* A}
    We have 
    \begin{equation*}
        D_{\text{cris}}(A_Z(b, x)) = x^* \mathcal{A}_Z, \hspace{0.5cm} \text{ for all } x \in X(\Q_p).
    \end{equation*}
\end{lemma}
Therefore, to compute the local $p$-adic height at $p$, $h_p(A_Z(b, x))$, it is sufficient to compute the Hodge filtration and Frobenius structure on this universal connection. Computing the Hodge filtration and Frobenius structure on $\mathcal{A}_Z$ is equivalent to computing two unipotent isomorphisms
\begin{equation*}
    \lambda^\star(x) : \Q_p \oplus V_{\dR} \oplus \Q_p(1) \xrightarrow \sim D_{\text{cris}}(A_Z(b, x)), \hspace{1cm} \text{ for } \star \in \{\text{\text{Fil}}, \phi\},
\end{equation*}
such that $\lambda^\phi$ respects the Frobenius action and $\lambda^{\text{\text{Fil}}}$ respects the Hodge filtration. Here, $\Q_p(1)$ denotes $D_{\text{cris}}(\Q_p(1))$, the image of $\Q_p(1)$ under $D_{\text{cris}}$, which is a one-dimensional filtered $\phi$-module with $\text{Fil}^{-1} = \Q_p$, $\text{Fil}^n = 0$ for all $n > -1$ and $\phi = 1/p$. With respect to a suitable basis for $D_{\text{cris}}(A_Z(b, x))$, the isomorphisms $\lambda^{\text{\text{Fil}}}, \lambda^\phi$ may be represented in $(2g + 2)$-block matrix form as
\begin{equation} \label{matrices: frob and fil}
\lambda^{\text{\text{Fil}}}(x) =
\begin{pmatrix}
        1 & 0 & 0 \\
        \boldsymbol{\alpha}_{\text{Fil}} & 1 & 0 \\
        \gamma_{\text{Fil}} & \boldsymbol{\beta}_{\text{Fil}}^\intercal  & 1
    \end{pmatrix}, \hspace{1cm} \lambda^\phi(x) = 
\begin{pmatrix}
        1 & 0 & 0 \\
        \boldsymbol{\alpha}_{\phi} & 1 & 0 \\
        \gamma_{\phi} & \boldsymbol{\beta}_{\phi}^\intercal  & 1
    \end{pmatrix}, 
\end{equation}
as in \cite{balakrishnan2019explicit} §4.5, §5.3 respectively. The isomorphism $\lambda^\phi$ is uniquely determined by the unique Frobenius structure on $\mathcal{A}_Z$ which gives a map $\phi^*\mathcal{A}_Z(b) \xrightarrow[]{\sim} \mathcal{A}_Z(b)$ sending $1$ to itself. In contrast, $\lambda^{\text{\text{Fil}}}$ is only well-defined up to the stabiliser of the Hodge filtration. Furthermore, the splitting $s$ of the Hodge filtration defines idempotents $s_1, s_2$ on $V_{\dR}$ with images $s(V_{\dR}/\text{Fil}^0V_{\dR})$ and $\text{Fil}^0V_{\dR}$ respectively. If our basis of $A_Z(b, x)$ has the property that $\boldsymbol{\alpha}_{\text{Fil}} = 0$, then \cite{balakrishnan2019explicit} Lemma 5.4 gives the following formula for the local $p$-adic height at $p$ 
\begin{equation} \label{lemma: h_p formula}
    h_p(A_Z(b,x)) = \chi_p(\gamma_\phi(x) - \gamma_{\text{Fil}}(x) - \boldsymbol{\beta}_{\phi}^\intercal \cdot s_1(\boldsymbol{\alpha}_\phi)(x) - \boldsymbol{\beta}_{\text{Fil}}^\intercal \cdot s_2(\boldsymbol{\alpha}_\phi)(x)).
\end{equation}

\par The quantities on the right hand side of (\ref{lemma: h_p formula}) depend only on the filtered $\phi$-module $D_{\text{cris}}(A_Z(b, x))$. Moreover, we have seen by the discussion after Lemma \ref{lemma: all local heights away from p are trivial} and Lemma \ref{lemma: eqn: D cris = x^* A}, that to compute the local height $h_p$, it suffices to construct the Hodge filtration and Frobenius structure of $x^*\mathcal{A}_Z$. Therefore, we would like to produce algorithms to compute all of the quantities on the right hand side of $(\ref{lemma: h_p formula})$ in terms of modular forms associated to $X$ without requiring an explicit plane model of the curve to produce a model-free quadratic Chabauty algorithm for $X$. The subject of \S \ref{sec: model-free theory} will be an algorithm to compute the matrix $\lambda^{\text{\text{Fil}}}$ in terms of modular forms, hence computing $\boldsymbol{\beta}_{\text{Fil}}$ and $\gamma_{\text{Fil}}$ in terms of modular forms. 

\subsection{Hodge filtration} \label{subsec: hodge}

\par We have just seen that the Hodge filtration is used to define the local $p$-adic height function $h_p$ as in (\ref{lemma: h_p formula}). In this section we will give an overview of the Hodge filtration calculation on $\mathcal{A}_Z$, as discussed in detail in \cite{balakrishnan2019explicit} §4. This is because our model-free algorithm, in \S \ref{sec: model-free theory}, is an adaptation of this algorithm and therefore closely follows the steps involved. We will explain that a key step of this process is computing a nice correspondence $Z$, and describe how to compute $Z$ in terms of a Hecke operator. This section will finish with Algorithm \ref{alg: hodge filtration} on how to compute the matrix $\lambda^{\text{Fil}}(x)$ as in (\ref{matrices: frob and fil}). Since this paper is concerned with the computation of the Hodge filtration on $\mathcal{A}_Z$, we omit details of the theory behind it and direct the reader to \S 4.3 of \cite{balakrishnan2019explicit}. 

\par Recall that $X$ is a smooth projective modular curve of genus $g > 1$ over $\Q$, such that $r= g$ and $\rho(J) > 1$. Let $Y \subset X$ be an affine open subset. Let $b$ be a rational point of $Y$. We follow \cite{balakrishnan2019explicit} \S 4 and specialise to the case that $Y = X \setminus \{\infty\}$, as this is what we always have in our examples. Therefore, we have $D \coloneqq X \setminus Y = \{\infty\}$, so $\#D = 1$, and $\{\infty\}$ is defined over $\Q$. Choose a set $\{\omega_0, ..., \omega_{2g - 1}\} \in \text{H}^0(Y, \Omega^1)$ of differentials on $Y$, so that they form a symplectic basis of $\text{H}_{\dR}^1(X)$, in the sense that the cup product is the standard symplectic form with respect to this basis. Let $\boldsymbol{\omega}$ denote the column vector $(\omega_0, ..., \omega_{2g-1})^\intercal$.

\par We recall that $Z$ is a nontrivial element in $\Ker(NS(J) \rightarrow NS(X))$, which is guaranteed to exist under our assumption that $\rho(J) > 1$, and we compute it as follows. For the algorithm, we require the action of $Z$ as a matrix on our basis of $\text{H}_{\dR}^1(X/\Q)$. We utilise the fact that $X$ is a modular curve and choose to construct the action of $Z$ in terms of the action of a certain Hecke operator, $T_q$, for $q$ a prime, with respect to our basis of $\text{H}_{\dR}^1(X/\Q)$. To compute the matrix of $T_q$, we first compute the matrix of Frobenius, $F$, via Tuitman's algorithm \cite{tuitman2016counting}, \cite{tuitman2017counting}, and then use the Eichler-Shimura relation:
\begin{equation*} 
    T_q = F + qF^{-1},
\end{equation*}
to obtain a matrix for the action of $T_q$ on $\text{H}_{\dR}^1(X/\Q)$. We choose the prime $q$ so that $T_q$ does not lie in $\Z \subset \text{End}(J)$, and in fact we can choose $q = p$. From here we choose to construct the matrix corresponding to the action of a nice correspondence $Z$ on $\text{H}_{\dR}^1(X/\Q)$ as a nontrivial polynomial in $T_p$, and we require that the corresponding matrix has trace equal to 0 so that $Z$ is guaranteed to be a nice correspondence. We will see in \S \ref{sec: model-free theory} that in the model-free algorithm, we first compute the action of $T_p$ on $q$-expansions and then use a change of basis to obtain the desired matrix for $Z$. This is because there is not a natural analogue of Tuitman's algorithm on $q$-expansions, so it is not straightforward to compute the matrix $F$.

\par We will now introduce the universal connection $\mathcal{A}_Z$. Firstly, universal properties give that the rank $2g + 2$ vector bundle $\mathcal{A}_Z$ has a connection, a Hodge filtration, and a Frobenius structure, as discussed in \cite{balakrishnan2019explicit} §4. We define $\mathcal{A}_Z$ to be a vector bundle with connection on $X$, $\nabla$, whose restriction to $Y$ is a trivial vector bundle with connection $\Lambda$ of the form
\begin{equation} \label{def of lambda}
   \Lambda \coloneqq - \begin{pmatrix}
        0 & 0 & 0 \\
        \boldsymbol{\omega} & 0 & 0 \\
        \eta & \boldsymbol{\omega}^\intercal  Z & 0
    \end{pmatrix}.
\end{equation}
The bottom left entry of $\Lambda$, $\eta$, is a differential of the 2nd kind on $X$, and we find that since $\#D = 1$, the differential $\eta$ can be taken to be zero. The connection $\nabla$ is uniquely defined in this way with respect to a trivialisation 
\begin{equation*}
    s_0 : \mathcal{O}_Y \otimes (\Q_p \oplus V_{\dR} \oplus \Q_p(1)) \xrightarrow{\sim} \mathcal{A}_Z \vert_Y, \hspace{1cm} \text{ so that we have } \hspace{0.2cm} s_0^{-1} \circ \nabla \circ s_0 = d + \Lambda,
\end{equation*}
and this trivialisation $s_0$ exists due a result of Kim \cite{kim2009unipotent} (see \cite{balakrishnan2019explicit} Theorem 4.1).

\par We now explain how to compute the Hodge filtration on $\mathcal{A}_Z$. The Hodge filtration on $\mathcal{A}_Z$ is uniquely determined by an isomorphism of filtered vector bundles
\begin{equation*}
    s^{\text{Fil}} : (\Q_p \oplus V_{\dR} \oplus \Q_p(1)) \otimes \mathcal{O}_Y \xrightarrow{\sim} \mathcal{A}_Z\vert_Y.
\end{equation*}
This morphism $s^\text{Fil}$ is uniquely determined by the vectors $\boldsymbol{\beta}_{\text{Fil}}$ and $\gamma_{\text{Fil}} \in \text{H}^0(Y, \mathcal{O}_Y)$ as entries of the matrix $\lambda^{\text{Fil}}$ in (\ref{matrices: frob and fil}) so that we in fact have $\lambda^{\text{Fil}} = s_0^{-1} \circ s^{\text{Fil}}$. Moreover, conditions imposed by Hadian's universal property \cite{hadian2011motivic} (\cite{balakrishnan2019explicit} Theorem 4.4) determine the matrix entries $\gamma_{\text{Fil}}$ and $\boldsymbol{\beta}_{\text{Fil}}$ uniquely as follows. At the point $\infty$, let 
\begin{equation*}
    s_\infty : \left((\Q_p \oplus V_{\dR} \oplus \Q_p(1)) \otimes \Q[[t]], d \right) \xrightarrow{\sim} \left( \mathcal{A}_Z\vert_{\Q[[t]]}, \nabla \right)
\end{equation*}
be a trivialisation of $\mathcal{A}_Z$ in a formal neighbourhood of $\infty$, with a local parameter $t$. Taking the composition of the trivialisation $s_0$ with the inverse of the trivialisation in the neighbourhood of $\infty$, $s_\infty$, we obtain a gauge transformation
\begin{equation*}
    G_\infty \coloneqq s_\infty^{-1} \circ s_0 \in \text{End}(\Q_p \oplus V_{\dR} \oplus \Q_p(1)) \otimes \Q((t)).
\end{equation*}
In addition, the gauge transformation, $G_\infty$, respects the filtration $\Lambda$, meaning that $G_\infty \circ (d + \Lambda) = d \circ G_\infty$, which is equivalent to $G_\infty$ satisfying the following equation
\begin{equation} \label{Gx and Lambda eqn}
    G_\infty^{-1}dG_\infty = \Lambda.
\end{equation}
Moreover, any such $G_\infty$ defines a trivialisation $s_\infty$. By recalling the definition of $\Lambda$ in (\ref{def of lambda}) and expanding out the equation in (\ref{Gx and Lambda eqn}) we see that $G_\infty$ is of the form
\begin{equation} \label{matrix: Gx def}
    G_\infty = \begin{pmatrix}
        1 & 0 & 0 \\
        \boldsymbol{\Omega} & 1 & 0 \\
        g_\infty & \boldsymbol{\Omega}^\intercal  Z & 1
    \end{pmatrix}, \hspace{1cm} \text{ where } \begin{cases}
        d\boldsymbol{\Omega} &= -\boldsymbol{\omega},  \\
        dg_\infty &= \boldsymbol{\Omega}^\intercal Z d\boldsymbol{\Omega}.
    \end{cases}
\end{equation}
Note that we have taken the definition of $G_\infty$ as in \cite{balakrishnan2023computation}, which is the minus of the definition taken in \cite{balakrishnan2019explicit}.

\par Finally, now that we have seen how to compute the matrix $G_\infty$, we can compute the matrix entries $\gamma_{\text{Fil}}, \boldsymbol{\beta}_{\text{Fil}}  \in \text{H}^0(Y, \mathcal{O}_Y)$. The vector $\boldsymbol{\beta}_{\text{Fil}}$ is in fact constant since Hadian's universal property \cite{hadian2011motivic} (\cite{balakrishnan2019explicit} Theorem 4.4) results in the condition that the functions in $\boldsymbol{\beta}_{\text{Fil}}$ extend to holomorphic functions on $X$, and so are constant. We can write  $\boldsymbol{\beta}_{\text{Fil}} \coloneqq (0_g, b_\text{Fil})^\intercal \in \Q_p^{2g}$, where $b_{\text{Fil}} = (b_g, ..., b_{2g-1})^\intercal \in \Q_p^g$, and let $N \coloneqq (0_g, 1_g)^\intercal \in M_{2g \times g}(\Q)$. Then, by \cite{balakrishnan2019explicit} \S 4, we have that $\gamma_{\text{Fil}}, \boldsymbol{\beta}_{\text{Fil}}$ must satisfy the following condition 
\begin{equation}  \label{hodge filtration equation gamma fil b fil}
    g_\infty + \gamma_{\text{Fil}} - b_{\text{Fil}}^\intercal N^\intercal \boldsymbol{\Omega}^\intercal - \boldsymbol{\Omega}^\intercal ZNN^\intercal \boldsymbol{\Omega} \in \Q[[t]].
\end{equation}
The existence and uniqueness of such $\gamma_{\text{Fil}}, \boldsymbol{\beta}_{\text{Fil}}$ follows from \cite{balakrishnan2019explicit}, Lemma 4.7. Moreover, we can determine the vector of constants $\boldsymbol{\beta}_{\text{Fil}}$ uniquely, and $\gamma_{\text{Fil}}$ is uniquely determined by the additional condition that $\gamma_{\text{Fil}}(b) = 0$. In summary, in this section we have given the following algorithm for computing the Hodge filtration on $\mathcal{A}_Z$, when $X = X_0^+(N)$, $N$ prime such that $g > 1$. 

\begin{algorithm} (Computing the Hodge filtration on $\mathcal{A}_Z$, Algorithm 3.12 \cite{balakrishnan2023quadratic}, Algorithm 5.20 \cite{balakrishnan2023computation}) \label{alg: hodge filtration}

Input: 
\begin{itemize}
    \item A prime $N$ such that $X_0^+(N)/\Q$ is a modular curve such that $g > 1$.
    \item A covering of $X$ by affine opens that are birational to a planar curve cut out by an equation that is monic in one variable, has $p$-integral coefficients and satisfies assumptions as found in Assumption 3.10 in \cite{balakrishnan2023quadratic}.
    \item A prime $p$ such that $T_p$ does not lie in $\Z \subset \End(J)$.
\end{itemize}
Output: the matrix $\lambda^{\text{Fil}}$.
\begin{enumerate}
    \item Compute an integral symplectic basis for $\text{H}^1_{\dR}(X/\Q)$.
    \item Compute a local coordinate $t$ at $\infty$.
    \item Expand $\boldsymbol{\omega}$ into a vector $\boldsymbol{\omega}$ of Laurent series.
    \item Compute the action of Frobenius on $\text{H}^1_{\dR}(X/\Q_p)$ using Tuitman's algorithm \cite{tuitman2016counting}, \cite{tuitman2017counting}. Use the Eichler-Shimura relation to compute the matrix of the action of the Hecke operator $T_p$ on $\text{H}^1_{\dR}(X/\Q_p)$.
    \item Compute the matrix of $Z \in  \Ker(NS(J) \rightarrow NS(X))$ acting on $\text{H}^1_{\dR}(X/\Q_p)$ in terms of $T_p$.
     \item Compute a vector $\boldsymbol{\Omega}$ that satisfies $d\boldsymbol{\Omega} = - \boldsymbol{\omega}$.
    \item Solve the system of equations for $g_\infty$ such that $dg_\infty = \boldsymbol{\Omega}^\intercal Z d\boldsymbol{\Omega}$.
    \item Compute the vector of constants $b_{\text{Fil}} = (b_g, ..., b_{2g-1}) \in \Q_p^g$ and the function $\gamma_{\text{Fil}}$ characterised by $\gamma_{\text{Fil}}(b) = 0$ and 
    \begin{equation*}
        g_\infty + \gamma_{\text{Fil}} - b_{\text{Fil}}^\intercal N^\intercal \boldsymbol{\Omega}^\intercal - \boldsymbol{\Omega}^\intercal ZNN^\intercal \boldsymbol{\Omega} \in \Q[[t]], 
    \end{equation*}
    where $N \coloneqq (0_g, 1_g)^\intercal \in M_{2g \times g}(\Q)$.
\end{enumerate}
    
\end{algorithm}

\begin{Remark} \label{remark: rho = g}
    We note that the inputs of Algorithm \ref{alg: hodge filtration} are not exactly how they are stated as in \cite{balakrishnan2023quadratic} Algorithm 3.12. By \cite{dogra2021quadratic} \S 4.1, we have that $\rho(J) = g$ when $J$ is the Jacobian of $X_0^+(N)$ for $N$ prime, and so the assumption that $g > 1$ gives that $\rho(J) > 1$ as required for the nice correspondence $Z$ to exist. Moreover, in \cite{balakrishnan2023quadratic}, it is stated that one computes all matrices of $Z_1, ..., Z_{g - 1}$ acting on $\text{H}^1_{\dR}(X/\Q_p)$, where $Z_1, ..., Z_{g - 1}$ is a basis of $\Ker(NS(J) \rightarrow NS(X))$, this is what requires the prime $p$ to be such that $T_p$ generates the Hecke algebra rather than $T_p$ does not lie in $\Z$. 
\end{Remark}



\section{Model-free Hodge filtration algorithm} \label{sec: model-free theory}

\par The focus of this section is the model-free algorithm for computing the Hodge filtration on $\mathcal{A}_Z$ for $X_0^+(N)$, the Atkin-Lehner quotient of the modular curve $X_0(N)$, for $N$ prime such that the genus, $g$, of $X$ is at least $2$. We saw in \S
\ref{subsec: hodge} that one of the first steps in the Hodge filtration, Algorithm \ref{alg: hodge filtration}, is to compute a basis of de Rham cohomology, $\text{H}^1_{\dR}(X/\Q)$, as the vector $\boldsymbol{\omega}$. In the model algorithm, this is defined in terms of the variables $x, y$ of the model of the curve. Therefore, the first step of a model-free algorithm for the Hodge filtration is to compute a basis of $\text{H}^1_{\dR}(X/\Q)$ in terms of $q$-expansions of modular forms, which will be dealt with in \S \ref{subsec: computing H1dR}. Then, in \S \ref{subsec: calculating hodge}, we will explain how to compute the nice correspondence $Z$, and the entries of the matrix $\lambda^{\text{Fil}}$ in terms of modular forms.

\par The holomorphic differential 1-forms on $X$ are a subspace of de Rham cohomology, and so we have $\text{H}^0(X_0^+(N), \Omega_{X_0^+(N)|\mathbb{C}}) \subset \text{H}^1_{\dR}(X/\Q)$. Therefore, since $\text{H}^0(X_0^+(N), \Omega_{X_0^+(N)|\mathbb{C}})$ has dimension $g$, and $\text{H}^1_{\dR}(X/\Q)$ has dimension $2g$, half of the basis of $\text{H}^1_{\dR}(X/\Q)$ comes from the holomorphic differential 1-forms. We obtain these elements as follows. Let $\mathcal{H}$ be the upper half plane, let $\tau \in \mathcal{H}$, and let $f$ be a weight $2$ cusp form on $X_0(N)$, so $f \in \mathcal{S}_2(\Gamma_0(N))$. Therefore, $f$ has a $q$-expansion of the form
\begin{equation*}
    f(\tau) = \sum\limits_{n = 1}^\infty a_n q^n, \hspace{0.5cm} q = e^{2\pi i \tau}, \tau \in \mathcal{H},
\end{equation*}
and from this the following $q$-expansion 
\begin{equation*}
    \omega_f \coloneqq f(\tau) \frac{dq}{q} = \sum\limits_{n = 1}^\infty a_n q^{n-1} dq,
\end{equation*}
defines an element of $\text{H}^0(X_0(N), \Omega_{X_0(N)|\mathbb{C}})$. Since $X_0^+(N)$ is a quotient of $X_0(N)$, we can think of $\text{H}^0(X_0^+(N), \Omega_{X_0^+(N)|\mathbb{C}})$ as a vector subspace of $\text{H}^0(X_0(N), \Omega_{X_0(N)|\mathbb{C}})$. Hence we get a corresponding subspace of $\mathcal{S}_2(\Gamma_0(N))$. This subspace is exactly the plus eigenspace of the action of the Atkin-Lehner involution on modular forms, and we denote this space as $\mathcal{S}_2^+(\Gamma_0(N))$.

\par Therefore, to compute a basis of $\text{H}^1_{\dR}(X/\Q)$, we first compute the space $\mathcal{S}_2^+(\Gamma_0(N))$, and then obtain a basis of $\text{H}^0(X_0^+(N), \Omega_{X_0^+(N)|\mathbb{C}})$ by sending a cusp form $f$ to $\omega_f \coloneqq f\frac{dq}{q}$. We denote the basis of $\text{H}^0(X_0^+(N), \Omega_{X_0^+(N)|\mathbb{C}})$ by $\{\omega_0, ..., \omega_{g-1} \}$, and the corresponding cusp forms by $\{ f_{\omega_0}, ..., f_{\omega_{g-1}}\}$, so that $\omega_i = f_{\omega_i}\frac{dq}{q}$ for $i \in \{0, ..., g-1\}$. The process of computing the basis $\{\omega_0, ..., \omega_{g-1} \}$ can be implemented in \texttt{Magma} \cite{MR1484478}. Therefore, it remains to construct the other half of the basis of $\text{H}^1_{\dR}(X/\Q)$. To do this, we will utilise the discriminant function, $\Delta$; this will be the focus of \S\ref{subsec: computing H1dR}. 

\par In \S \ref{subsec: computing H1dR}, we will show that our method always produces a basis of $\text{H}^1_{\dR}(X/\Q)$ and we will verify this by computing the cup product of pairs of our basis elements to compute a cup product matrix. If our cup product matrix has nonzero determinant, then we have successfully constructed a basis of $\text{H}^1_{\dR}(X/\Q)$ in terms of $q$-expansions of modular forms. We will use Serre's cup product formula:
\begin{align*}
    \text{H}^1_{\dR}(X/\Q_p) \times \text{H}^1_{\dR}(X/\Q_p) &\rightarrow \Q_p \\ 
    ([\mu_1], [\mu_2]) &\mapsto [ \mu_1 \cup \mu_2] = \sum_{Q \in X(\mathbb{C}_p)} \text{Res}_Q \left( \mu_2 \int \mu_1 \right).
\end{align*}    

\par Therefore, if we calculate the cup product matrix of our basis and find that its determinant is nonzero, then we have successfully found a basis of $\text{H}^1_{\dR}(X/\Q)$. Since the differentials $\{\omega_0, ..., \omega_{g-1}\}$ are holomorphic, by the definition of the cup product formula we will have $[\omega_i \cup \omega_j] = 0$ for all $i, j, \in \{0, ..., g-1\}$.
Previous successful implementations of quadratic Chabauty in the literature (see for example \cite{balakrishnan2023quadratic}) rely on having a symplectic basis as an input as this makes some later calculations simpler. Therefore we will use change of basis to make our basis a symplectic one, where symplectic simply means that its cup product matrix is equal to 
\begin{equation} \label{def: symp cup product matrix}
    C \coloneqq \begin{pmatrix}
        0_g & I_g \\ -I_g & 0_g
    \end{pmatrix}, \hspace{0.5cm} \text{where }  \hspace{0.5cm} C_{i,j} = [\omega_i] \cup [\omega_j], \text{ for } i, j \in \{0, ..., 2g-1\}.
\end{equation}

\par Once we have a symplectic basis of $\text{H}^1_{\dR}(X/\Q)$, we will then need to compute the action of the Hecke operator $T_p$ on $q$-expansions of differentials, in order to compute the matrices $Z$, and $\Lambda$ in terms of this basis of $\text{H}^1_{\dR}(X/\Q)$ as $q$-expansions. We will then explain how to compute the values $\boldsymbol{\beta}_{\text{Fil}}, \gamma_{\text{Fil}}$. This will be the focus of \S\ref{subsec: calculating hodge}.

\subsection{Computing a symplectic basis of $\text{H}^1_{\dR}(X/\Q)$} \label{subsec: computing H1dR}

\par We have already seen how we obtain the first half of a basis for $\text{H}^1_{\dR}(X/\Q)$ from $\mathcal{S}_2^+(\Gamma_0(N))$, and that we will use Serre's cup product formula to verify when we have constructed a basis of $\text{H}^1_{\dR}(X/\Q)$. Since Serre's cup product formula sums over the poles of the differentials, we would like to construct differentials that we know only have a pole at the unique cusp $\infty$ of $X_0^+(N)$. To do this, we will construct a function $f_{\dR}$ with poles only at infinity and then multiply this by the holomorphic differentials. We will then compute the cup product matrix of the set $\{\omega_0, ..., \omega_{g-1}, f_{\dR} \cdot \omega_0, ..., f_{\dR} \cdot \omega_{g-1}\}$ to verify that this set is a basis for $\text{H}^1_{\dR}(X/\Q)$. We will construct $f_{\dR}$ as a quotient of two weight 12 cusp forms on $X_0^+(N)$.  In this section we give a suitable construction of the function $f_{\dR}$ for $X_0^+(N)$, and prove that our method always produces a basis. We begin by giving a function with zeroes only at $\infty$, using the discriminant function $\Delta \in \mathcal{S}_{12}(\text{SL}_2(\mathbb{Z}))$.

\begin{prop} \label{prop: properties of sqrt delta w_N(delta)}
    For any odd prime $N$, the function $\sqrt{\Delta \cdot w_N(\Delta)}$ is a weight 12 cusp form for $\Gamma_0(N)$ in the plus eigenspace for the Atkin-Lehner involution, i.e the function $\sqrt{\Delta \cdot w_N(\Delta)}$ is in $\mathcal{S}_{12}^+(\Gamma_0(N))$. Moreover, the function $\sqrt{\Delta \cdot w_N(\Delta)}$ only vanishes at the unique cusp at $\infty$ on $X_0^+(N)$.

\begin{proof}
    By definition, $\Delta = \eta(z)^{24}$, where $\eta(z) = q^{1/24}\prod\limits_{n=1}^\infty (1 - q^n)$ is the Dedekind eta-function. We also have $w_N(\Delta(z)) = \Delta(w_n\cdot z) = \Delta(-\frac{1}{Nz}) = \Delta(Nz) = \eta(Nz)^{24}$, where the first two equalities follow from the definition of $w_N$, the third is because $\Delta$ is a cusp form of weight 12, and the final equality follows from the definition of $\Delta$. Therefore, we can rewrite $\sqrt{\Delta \cdot w_N(\Delta)}$ as $\eta(z)^{12}\eta(Nz)^{12}$. By  \cite{rouse2015spaces} (2), a function $f = \prod_{\delta \vert N} \eta(\delta z)^{r_\delta}$ is a modular form for $\Gamma_0(N)$ of weight $k = \frac{1}{2}\sum_{\delta \vert N} r_{\delta}$, if we have 
    \begin{align*}
        &\sum\limits_{\delta \vert N} \delta r_\delta \equiv 0 \text{ mod } 24, \hspace{0.4cm} \sum\limits_{\delta \vert N} \frac{N}{\delta} r_\delta \equiv 0 \text{ mod } 24, \text{ and } \\
        &\prod\limits_{\delta \vert N} \delta^{r_{\delta}} \text{ is the square of a rational number}.
    \end{align*}
    Therefore, since the function $\eta(z)^{12}\eta(Nz)^{12}$ has pairs of values of $(\delta, r_\delta)$ equal to $ \{(1, 12), (N, 12)\}$, all of the above conditions are satisfied and so $\eta(z)^{12}\eta(Nz)^{12}$ is a modular form for $\Gamma_0(N)$ of weight 12. Finally, $\eta(z)^{12}\eta(Nz)^{12}$ is a cusp form and invariant under the action of $w_N$ and hence $\sqrt{\Delta \cdot w_N(\Delta)} \in \mathcal{S}_{12}^+(\Gamma_0(N))$. The cusp form $\Delta$ is nonzero for all $\tau \in \mathcal{H}$ (see, for example Corollary 1.4.2 \cite{diamond2005first}), therefore $w_N(\Delta)$ vanishes only at $\infty$, and so the product $\sqrt{\Delta \cdot w_N(\Delta)}$ has a zero only at $\infty$.
\end{proof}
\end{prop}

\par We will use the function $\sqrt{\Delta \cdot w_N(\Delta)}$ to define the function $f_{\dR}$ as follows. Write $val(f)$ for the valuation of a $q$-expansion $f$ as an element of $\Q((q))$. We let $d \coloneqq \text{dim } \mathcal{S}_{12}^+(\Gamma_0(N))$ and write $\{s_1, ..., s_d\}$ to be a basis for $\mathcal{S}_{12}^+(\Gamma_0(N))$ such that $val(s_i) < val(s_{i + 1})$ for all $i \in \{1, ..., d-1\}$. We can compute such a basis on \texttt{Magma} \cite{MR1484478} and it gives $s_d = a_\Delta\sqrt{\Delta \cdot w_N(\Delta)}$, where $a_{\Delta}$ is a constant so that the leading term of $s_d$ is 1. We have $val(\sqrt{\Delta \cdot w_N(\Delta)}) = \frac{1}{2} + \frac{N}{2}$, and so $val(s_d) = \frac{N+1}{2}$. We will construct $f_{\dR}$ by choosing a weight 12 cusp form in the plus eigenspace for the Atkin-Lehner involution, $s_{d - j_{\dR}}$, for a positive integer $j_{\dR} \in \{1, ..., d-1\}$ satisfying some extra conditions, and divide it by the function $a_\Delta \sqrt{\Delta \cdot w_N(\Delta)}$. Therefore, the function $f_{\dR}$ will have a pole only at $\infty$ due to Proposition \ref{prop: properties of sqrt delta w_N(delta)}.

\par We now discuss how to choose the integer $j_{\dR} \in \mathbb{N}$. We want to choose $j_{\dR} \in \{1, ..., d-1\}$ so that taking $f_{\dR} = \frac{s_{d - j_{\dR}}}{s_d}$ results in the set $\{\omega_0, ..., \omega_{g-1}, f_{\dR}\cdot\omega_0, ..., f_{\dR}\cdot\omega_{g-1}\}$ being a basis of $\text{H}^1_{\dR}(X/\Q)$. For this, we require, for all $i \in \{0, ..., g-1\}$, that 
\begin{equation} \label{eqn: condition on j}
    val\left(\frac{s_{d - j_{\dR}}}{s_d}f_{\omega_i}\right) \leq -1.
\end{equation}
This is because if $val\left(\frac{s_{d - j_{\dR}}}{s_d}f_{\omega_i}\right) \geq 1$ for some $i \in \{0, ..., g-1\}$, then $\frac{s_{d - j_{\dR}}}{s_d}f_{\omega_i} \in \mathcal{S}_2^+(\Gamma_0(N))$ and so under the map to differentials, it produces an element of $\text{H}^0(X, \Omega^1)$. Therefore, the set $\{\omega_0, ..., \omega_{g-1}, f_{\dR}\cdot\omega_0, ..., f_{\dR}\cdot\omega_{g-1}\}$ could not give a basis of $\text{H}^1_{\dR}(X/\Q)$, as the differential $\frac{s_{d - j_{\dR}}}{s_d}f_{\omega_i}\frac{dq}{q}$ would lie in the span of $\{\omega_0, ..., \omega_{g-1}\}$. Moreover, if $val\left(\frac{s_{d - j_{\dR}}}{s_d}f_{\omega_i}\right) = 0$ for some $i \in \{0, ..., g-1\}$, then the computation of the cup product formula will fail since we will have to integrate a nonzero residue term. Therefore, we must choose $j_{\dR} \in \mathbb{N}$ such that (\ref{eqn: condition on j}) holds. 

\begin{prop} \label{prop: there exists a j}
    For any prime number $N$, there exists an integer $j \in \{1, ..., d-1\}$ such that for all $i \in \{0, ..., g-1\}$ we have
    \begin{equation*}
        val\left(\frac{s_{d - j}}{s_d}f_{\omega_i}\right) \leq -1.
    \end{equation*}

\begin{proof}
We have $val(s_d) = \frac{N+1}{2}$. We defined the set $\{\omega_0, ..., \omega_{g-1}\}$ to be a basis of $\text{H}^0(X, \Omega^1)$ ordered by increasing valuation. Therefore, for all $j \in \{1, ..., d-1\}$, we have the following upper bound for all $i \in \{0, ..., g-1\}$. 
    \begin{align*}
    val\left(\frac{s_{d - j}}{s_d}f_{\omega_i}\right) &\leq val(s_{d - j}) - val({s_d}) + \max_{i=0, ..., g-1} (val(f_{\omega_i})) \\
    &= val(s_{d - j}) - \frac{N+1}{2} + val(f_{\omega_{g-1}}).
\end{align*}

\par Take $j = d-1$, then $val(s_{d - j}) = val(s_1) = 1$, because otherwise we will have $val(f) \geq 2$ for all $f \in \mathcal{S}_{12}^+(\Gamma_0(N))$ which gives a contradiction since any line bundle of degree at least $2g$ on  curve of genus $g$ is basepoint free. Therefore, we have $val(s_1) = 1$. All that is left is to show that 
\begin{equation} \label{eqn: ineq for valuation of f}
val(f_{\omega_{g-1}}) \leq \frac{N-3}{2}.
\end{equation}
\par We will prove this by using the valence formula for congruence subgroups to prove that (\ref{eqn: ineq for valuation of f}) holds (see, for example, \cite{masdeu2015modular} Theorem 2.6.1). Since $f \in \mathcal{S}_2^+(\Gamma_0(N))$, $f$ lies in the larger space $ \mathcal{S}_2(\Gamma_0(N))$. Let $\text{PSL}_2(\mathbb{Z}) = \text{SL}_2(\mathbb{Z}) / \{\pm 1\}$, and $\overline{\Gamma}$ is the image of $\Gamma$ in $\text{PSL}_2(\mathbb{Z})$, where $\Gamma$ is a congruence subgroup. Also, define $v_z(f)$ as the valuation of $f$ at the point $z \in \mathcal{H} \cup \{ \infty \}$, $n_\Gamma(z)$ is the size of the stabiliser of $z$ in $\overline{\Gamma}$, $v_{P, \Gamma}(f)$ is the valuation of $f$ at the cusp $P$, and $Cusps(\Gamma)$ are the set of cusps of $\Gamma$. The valence formula for congruence subgroups gives that for a nonzero meromorphic function on $\mathcal{H} \cup \{\infty\}$ which is weakly modular of weight $k$ for the congruence subgroup $\Gamma$, we have the following
\begin{equation} \label{eqn: valence formula}
    \sum\limits_{z \in \Gamma\backslash \mathcal{H}} \frac{v_z(f)}{n_{\Gamma}(z)} + \sum\limits_{P \in Cusps(\Gamma)} v_{P, \Gamma}(f) = \frac{k}{12}[\text{PSL}_2(\mathbb{Z}): \overline{\Gamma}]. 
\end{equation}

\par Since $f$ is a cusp form, we have $v_z(f) \geq 0$ for all $z \in \Gamma\backslash \mathcal{H}$, and writing $w_N(\infty)$ for the second cusp of $X_0(N)$, we have $v_{w_N(\infty), \Gamma_0(N)}(f) > 0$. We also have $[\text{PSL}_2(\mathbb{Z}): \overline{\Gamma_0(N)}] = [\text{SL}_2(\mathbb{Z}): \Gamma_0(N)] = N + 1$. Therefore, since $val(f_{\omega_{g-1}}) = v_{\infty, \Gamma_0(N)}(f)$, the valence formula (\ref{eqn: valence formula}) gives
\begin{equation*}
    val(f_{\omega_{g-1}}) = \frac{N+1}{6} - \left(\sum\limits_{z \in \Gamma\backslash \mathcal{H}} \frac{v_z(f)}{n_{\Gamma}(z)} + v_{w_N(\infty), \Gamma_0(N)}(f) \right) \leq \frac{N+1}{6}.
\end{equation*}
Finally, we have $\frac{N+1}{6} \leq \frac{N-3}{2}$ if and only if $N \geq 5$, which we have since we have assumed that the genus of $X_0^+(N)$ is at least 2, and so (\ref{eqn: ineq for valuation of f}) holds.
\end{proof}
\end{prop}

\begin{Remark}
    A different approach would be to set $j = 1$ and then take $m$, a positive integer, large enough such that $(\frac{s_{d-1}}{s_d})^mf_{\omega_i}$ has valuation less than or equal to $-1$ for all $i \in \{0, ..., g-1\}$.
\end{Remark}

\par Therefore, by Proposition \ref{prop: there exists a j}, we know that there always exists a $j \in \{1, ..., d-1\}$ such that (\ref{eqn: condition on j}) holds. Our strategy will be to define $j_{\dR}$ to be the smallest such $j \in \mathbb{N}$. Once we have found such a $j_{\dR} \in \mathbb{N}$ we define $f_{\dR}$ to be 
\begin{equation*}
    f_{\dR} = \frac{s_{d - j_{\dR}}}{s_d} = \frac{s_{d - j_{\dR}}}{ a_{\Delta}\sqrt{\Delta \cdot w_N(\Delta)}}.
\end{equation*}
We then compute the cup product matrix for the set $\{\omega_0, ..., \omega_{g-1}, f_{\dR}\cdot\omega_0, ..., f_{\dR}\cdot\omega_{g-1}\}$ to verify that we have indeed found a basis.

\par We end this section by proving that this method of constructing a basis of $q$-expansions for $\text{H}^1_{\dR}(X/\Q)$ will always work. We begin with the following two results, which show that any differential of the 2nd kind on $X$ can be written in $\text{H}^1_{\dR}(X/\Q)$ as a product of a global differential 1-form and a function with poles only at infinity of bounded order.

\begin{prop} \label{prop: there exists a g ord(w0) =< ord(w) + m infty}
    Let $\omega$ be a differential of the 2nd kind on $X$, let $\omega_0 \in \text{H}^0(X, \Omega)$ be such that $ord_\infty(\omega_0) = 0$, and let $div(\omega_0) = \sum\limits_{i=1}^\ell n_i P_i$ be the divisor associated to $\omega_0$ such that $n_i > 0$ for all $i \in \{0, ..., \ell\}$. Let $m \in \textbf{N}$ be such that $m > 4g - 4 + \ell$. Then, there exists a function $g \in \C(X)$ such that $\omega - dg \in \text{H}^0(X, \Omega(m\infty))$ and $div(\omega_0) \leq div(\omega - dg) + m \infty$.

\begin{proof}
    \par For any differential of the 2nd kind on $X$, $\omega$, there exists a meromorphic function $f$ such that $\omega - df$ has a pole at $\infty$ only. Therefore, without loss of generality, in $\text{H}^1_{\dR}(X/\Q)$ we assume that $\omega$ has a pole at $\infty$ only, and we can bound the size of this pole at $\infty$ as follows. We claim that we can choose a function $g_1 \in \C(X)$ such that $\omega - dg_1$ has a pole at $\infty$ only of order at most $2g$. Let $m_0 > 0, m_0 \in \mathbf{N}$. Suppose 
    \begin{equation} \label{eqn: dimension increase}
        \dim \text{H}^0(X, \mathcal{O}(m_0\infty)) = \dim \text{H}^0(X, \mathcal{O}((m_0-1)\infty) + 1. 
    \end{equation}
    Then, for every differential $\omega \in \text{H}^0(X, \Omega((m_0+1)\infty))$, there exists a function $g \in \text{H}^0(X, \mathcal{O}(m_0\infty))$ such that $\omega - dg \in \text{H}^0(X, \Omega(m_0\infty))$. Indeed, if the dimension increases, then there exists a function $g_0$ with a pole only at $\infty$ of order exactly $m_0$. Then, by taking its derivative and rescaling we can remove the leading term of the Taylor expansion of $\omega$ at $\infty$ by subtracting $dg$, where $g = \lambda g_0$, where $\lambda$ is the suitable scalar so that the leading term of $\omega$ is removed. Moreover, by Riemann-Roch, if $m_0 > 2g - 1$, then (\ref{eqn: dimension increase}) holds, and since we have assumed that $m > 4g - 4 + \ell$ we have $m > 2g - 1$. As a result, if we inductively apply the argument above, then if $\omega$ only as a pole at $\infty$, there exists a $g_1$ with only a pole at $\infty$ such that $\omega - dg_1$ only has a pole at $\infty$ of order at most $2g$. This means that $ord_\infty(\omega_0) \leq ord_\infty(\omega - dg_1) + m\infty$. Therefore, we have shown that, without loss of generality, in $\text{H}^1_{\dR}(X/\Q)$ we have that $\omega$ has a pole only at $\infty$ of order at most $2g$. Additionally, we have assumed that $P_i \neq \infty$ for all $i \in \{1, ..., \ell\}$, and so it remains to prove that for all $i \in \{1, ..., \ell\}$, we can construct a function $g_2 \in \C(X)$ such that $n_i = ord_{P_i}(\omega_0) \leq ord_{P_i}(\omega - dg_2)$, assuming that $m > 4g - 4 + \ell$. 

    \par Let $t_i$ be a parameter at $P_i$ for all $i \in \{1, ..., \ell\}$. Let the Taylor expansion of $\omega$ at $P_i$ be $ \sum\limits_{j=0}^\infty a_{i,j} t_i ^j dt_i$, and let the Taylor expansion of a function $g$ at $P_i$ be written as $\sum\limits_{j=0}^\infty b_{i,j}(g) t_i ^j$, so that the coefficients, $b_{i, j}(g)$, depend on the function $g$. We want to prove that we can always find a function $g_2$ such that 
    \begin{equation} \label{eqn: w and dg coefficients}
        j \cdot b_{i,j}(g_2) = a_{i,j-1} \hspace{0.4cm} \text{ for all }i \in \{1, ..., \ell\} \text{ and for all } j \in \{1, ..., n_i\}.
    \end{equation}
    This is because if (\ref{eqn: w and dg coefficients}) holds, then $\omega - dg_2$ has $\sum\limits_{j=n_i}^\infty (a_{i, j} - (j+1)\cdot b_{i, j+1}(g_2))t_i^j dt_i$ as a Taylor expansion at $P_i$ for all $i \in \{1, ..., \ell\}$. Therefore, $\omega - dg_2$ will vanish at $t_i$ of order at least $n_i$, and so $n_i \leq ord_{P_i}(\omega - dg_2)$ for all $i \in \{1, ..., \ell\}$. To show that we can always find such a function $g_2$, we first define a map $\delta$ 
    \begin{align*}
        \delta : \text{H}^0(X, \mathcal{O}(m\infty)) &\longrightarrow \C ^{2g - 2 + \ell},  \\
        g &\mapsto ((b_{ij}(g))_{j=0}^{n_i})_{i=1}^\ell,
    \end{align*}
    where we have used the fact that $\sum\limits_{i=1}^\ell (n_i + 1) = 2g - 2 + \ell$ because $div(\omega_0)$ is a canonical divisor. It is enough to prove that the map $\delta$ is surjective, so that we can always find a function $g_2$ such that (\ref{eqn: w and dg coefficients}) holds. We have $\Ker(\delta) = \text{H}^0(X, \mathcal{O}(m\infty - \sum\limits_{i=1}^\ell (n_i + 1)P_i))$, which are functions with poles only at $\infty$ of order at most $m$, and zeroes at $P_i$ of order at least $n_i + 1$ for all $i \in \{1, ..., \ell\}$. For $\delta$ to be surjective, it is sufficient to prove that
    \begin{equation} \label{eqn: condition for delta to be surjective}
       \dim\text{H}^0(X, \mathcal{O}(m\infty)) =\dim\text{H}^0(X, \mathcal{O}(m\infty - \sum\limits_{i=1}^\ell (n_i + 1)P_i)) + 2g - 2 + \ell.
    \end{equation}
    Let $D_1 \coloneqq m\infty - \sum\limits_{i=1}^\ell (n_i + 1)P_i$. Since we have assumed that $m > 4g - 4 + \ell$, we have $deg(D_1) > 2g - 2$, and therefore $\dim\text{H}^0(X, \mathcal{O}(D_1)) = deg(D_1) - g + 1$, by Riemann-Roch. Let $D_2 \coloneqq \sum\limits_{i=1}^\ell (n_i + 1)P_i$. Then $D_1 + D_2 = m\infty$, and since $deg(D_2) > 0$, we similarly have that $\dim \text{H}^0(X, \mathcal{O}(D_1 + D_2)) = deg(D_1 + D_2) + g - 1$. Therefore, taking the difference of these dimension formulas gives 
    \begin{equation*}
       \dim\text{H}^0(X, \mathcal{O}(D_1 + D_2)) - \dim \text{H}^0(X, \mathcal{O}(D_1)) = deg(D_2) = 2g - 2 + \ell. 
    \end{equation*}
    Therefore, we have the equality in (\ref{eqn: condition for delta to be surjective}) and so the map $\delta$ is surjective. Taking $g \coloneqq g_1 + g_2$ gives the function we require.
\end{proof}
\end{prop}

\par The condition in Proposition \ref{prop: there exists a g ord(w0) =< ord(w) + m infty} that for a point $P \neq \infty$ on $X$ we have $ord_P(\omega_0) \leq ord_P(\omega - dg)$ guarantees that the quotient $\frac{\omega - dg}{\omega_0}$ has no pole at any point $P \neq \infty$. In addition, the assumption that $ord_\infty(\omega_0) = 0$ gives that the quotient $\frac{\omega - dg}{\omega_0}$ has a pole only at $\infty$ of order at most $m$. This leads to the following important result.

\begin{corollary} \label{cor: w of 2nd kind can be written f w_0}
    Let $\omega, \omega_0$ and $m \in \mathbb{N}$ be as in Proposition \ref{prop: there exists a g ord(w0) =< ord(w) + m infty}. Then there exists a function $f \in \text{H}^0(X, \mathcal{O}(m\infty))$ and $\omega_0 \in \text{H}^0(X, \Omega^1)$ such that $[\omega] = [f \cdot \omega_0]$ in $\text{H}^1_{\dR}(X/\Q)$.

\begin{proof}
By Proposition \ref{prop: there exists a g ord(w0) =< ord(w) + m infty}, there exists a function $g \in \C(X)$ such that $\omega - dg$ only has a pole at $\infty$ of order $m$, and $div(\omega_0) \leq div(\omega - dg) + m\cdot [\infty]$. This means that there exists a function $f \in \text{H}^0(X, \mathcal{O}(m\infty))$ such that $\omega - dg = f \cdot \omega_0$. As a result, $[\omega] = [f\cdot \omega_0]$ in de Rham cohomology, as required.
\end{proof}
\end{corollary}

\par Corollary \ref{cor: w of 2nd kind can be written f w_0} gives a method to construct a basis of $\text{H}^1_{\dR}(X/\Q)$ by constructing functions with poles only at infinity of sufficiently large order. In this section, we have done this by taking a quotient of weight 12 cusp forms of the form $\frac{s_{d-j}}{a_{\Delta}\sqrt{\Delta \cdot w_N(\Delta)}}$ for $j \in \{1, ..., d-1\}$. In particular, this allows us to produce functions with poles only at infinity of order at most $-val(\frac{s_1}{a_{\Delta}\sqrt{\Delta \cdot w_N(\Delta)}}) = \frac{N-1}{2}$. We want to show that $\frac{N-1}{2} > 4g - 4 + \ell$, so that this construction satisfies the conditions of Proposition \ref{prop: there exists a g ord(w0) =< ord(w) + m infty}, and since $\ell \leq 2g - 2$, it is enough to show that $\frac{N-1}{2} > 6g - 6$. We claim that
\begin{equation} \label{eqn: (N+1)/2 = d + g}
    \frac{N+1}{2} = d + g.
\end{equation}
Indeed, let $g_0(N)$ be the genus of $X_0(N)$, so $g \leq g_0(N)$, and let $D$ be the divisor $D = \frac{N+1}{2}[\infty]$. Then, by Riemann-Hurwitz (\cite{diamond2005first} Theorem 3.1.1), we have $g_0(N) \leq \frac{N+1}{12}$, meaning that $2g - 2 < \frac{N+1}{2}$, and so by Riemann-Roch we have $\ell(D) = \frac{N+1}{2} - g + 1$. Additionally, using Riemann-Roch to compute the space of cusp forms, (such as in (3.12) of \cite{diamond2005first}), we have $d = \ell(D) - \epsilon_\infty$, where $\epsilon_\infty$ is the number of cusps for $X_0^+(N)$. Therefore, since $\epsilon_\infty = 1$ as $N$ is prime, we obtain  $\frac{N+1}{2} = d + g$. Following the equality in (\ref{eqn: (N+1)/2 = d + g}), we therefore want to prove that $d > 5g - 5$. By \cite{diamond2005first} Theorem 3.5.1, we have $\dim \mathcal{S}_{12}(\Gamma_0(N)) \geq 11g_0(N) - 1$. Moreover, by \cite{martin2018refined} Theorem 2.2, we have that we have $\dim \mathcal{S}_{12}^{new, +}(\Gamma_0(N)) \geq \frac{1}{2} \dim \mathcal{S}_{12}^{new}(\Gamma_0(N))$, and so $d \geq \frac{1}{2}\dim \mathcal{S}_{12}(\Gamma_0(N))$, as $N$ is prime. We therefore have $d > 5g - 5$ as required. This means that our method of taking quotients of the form $\frac{s_{d-j}}{a_{\Delta}\sqrt{\Delta \cdot w_N(\Delta)}}$ for $j \in \{1, ..., d-1\}$ does indeed produce functions with poles only at infinity of sufficiently large order. 

\par Write $n \coloneqq \frac{N - 1}{2}$, for the maximum pole order we can construct from our basis of $\mathcal{S}_{12}^+(\Gamma_0(N))$. Since $n > 2g-2$, by the equality in (\ref{eqn: (N+1)/2 = d + g}) and $d > 5g - 5$, by a similar argument to the proof of Proposition \ref{prop: there exists a g ord(w0) =< ord(w) + m infty}, we have that $\dim\text{H}^0(X, \mathcal{O}(n\infty)) = n - g + 1 = d$. Therefore, if we construct a basis $f_1, ..., f_d$ of $\text{H}^0(X, \mathcal{O}(n\infty))$ such that $ord_\infty(f_i) < ord_\infty(f_{i+1})$ for all $i \in \{1, ..., d-1\}$, where we write $a_i \coloneqq ord_\infty(f_i)$ for all $i \in \{1, ..., d\}$. Given such a basis of $\text{H}^0(X, \mathcal{O}(n\infty))$, we can construct a basis of $\text{H}^0(X, \mathcal{O}(2n\infty))$ by taking products of our basis $f_1, ..., f_d$ by the following. First, we have $\{2g, ..., n\} \subset \{a_1, ..., a_d\}$. It is sufficient to prove that every $\tilde{n} \in \{n + 1, ..., 2n\}$ there exists $a, b \in \{2g, ..., n\}$ such that $\tilde{n} = a + b$. This is equivalent to requiring that $4g \leq n + 1$, which can be written as $4g \leq \frac{N+1}{2} = d + g$, which we have since $d \geq 11 g_0(N) + \frac{1}{2}(g_0(N) - 1)$. Therefore, once we have a basis of $\text{H}^0(X, \mathcal{O}(n\infty))$ we can construct a basis of $\text{H}^0(X, \mathcal{O}(2n\infty))$, and can use this to build differentials with poles only at infinity of any order.

\par By the discussion above, our method of computing a basis of $\text{H}^1_{\dR}(X/\Q)$ in terms of multiplying global differentials by functions with poles only at infinity will always work, if we take $f_{\dR} = \frac{s_1}{s_d}$. Moreover, by computing the cup product, we have found that taking $j_{\dR}$ as in (\ref{eqn: condition on j}) always produces a basis in our examples in \S \ref{sec: example}. We will perform a change of basis on our current basis of $\text{H}^1_{\dR}(X/\Q)$, denoted by $\{\omega_0, ..., \omega_{g-1}, \eta_g, ..., \eta_{2g-1}\}$, (so that $\eta_{i+g} \coloneqq f_{\dR}\cdot \omega_i$ for $i=0, ..., g-1$), to construct a new basis $\{\omega_0, ...,\omega_{g-1}, \omega_g, ..., \omega_{2g-1}\}$ such that its cup product matrix is equal to the matrix $C$, as defined in (\ref{def: symp cup product matrix}). We write $\boldsymbol{\omega} = (\omega_0, ..., \omega_{2g-1})$ for the basis written as entries of a vector. We have written code to implement such a change of basis in \texttt{Magma} \cite{MR1484478}.

\subsection{Calculating the Hodge filtration} \label{subsec: calculating hodge}

\par One of the main steps in computing the Hodge filtration and the Frobenius structure is to compute a certain unipotent vector bundle with connection. If $Y \subset X$ is an affine open with coordinate ring $R = \mathcal{O}(Y)$, then a unipotent vector bundle with connection is uniquely determined by the data of a $(2g+2) \times (2g+2)$ matrix $\Lambda$ with entries in $\Omega_{R|K}$, as first defined in (\ref{def of lambda}). To compute the entries of $\Lambda$, we need to compute the action of the Hecke Operator $T_p$ on $\text{H}^1_{\dR}(X/\Q)$, and the nice correspondence $Z$.

\par We start by computing the action of $T_p$ on $\text{H}^1_{\dR}(X/\Q)$ as a $2g \times 2g$ matrix. The action of $T_p$ preserves the subspace $\text{H}^0(X, \Omega^1) \subset \text{H}^1_{\dR}(X/\Q)$. Therefore in block form it will look like
\begin{equation*}
    T_p = \begin{pmatrix}
        * & * \\ 0_{g} & *
    \end{pmatrix}.
\end{equation*}
Denote the entries of the matrix $T_p$ as $(t_{j,i})_{ij}$ for $i, j \in \{ 0, ..., 2g-1\}$. We now explain how to compute the remaining values $t_{j, i}$. For each $i \in \{0, ..., 2g-1\}$ we would like to compute constants  $t_{i, j}$ so that 
\begin{equation*}
    T_p(\omega_i) = \sum\limits_{j = 0}^{2g-1} t_{i, j} \omega_j,
\end{equation*}
so that the action of $T_p$ on a given $\omega_i$ for $i \in \{0, ..., 2g-1\}$ is written as a linear combination of our symplectic basis of $\text{H}^1_{\dR}(X/\Q)$. We will compute these values $t_{i, j}$ as follows. Let $M \coloneqq ([T_p \omega_i] \cup [\omega_j])_{ij}$ for $i, j, \in \{0, ..., 2g-1\}$ be the cup product matrix between the symplectic basis and its image under $T_p$. We will first compute the matrix $M$ by computing the action of $T_p$ on $q$-expansions of differentials as in (\ref{eqn: Tp induced map}), and then use the following result.

\begin{prop} \label{prop: Tp = (MC^-1)^T}
    We have $T_p = (MC^{-1})^\intercal $.

\begin{proof}
    We have $C^{-1} = -C = C^\intercal$. Therefore, from the definitions of $M, C$ we have
    \begin{align*}
        MC^{-1} &= \begin{pmatrix}
            [T_p\omega_0 \cup \omega_0] & ... &[T_p\omega_0 \cup \omega_{2g-1}] \\
            \vdots & \ddots & \vdots \\
            [T_p\omega_{2g-1} \cup \omega_0] & ... & [T_p\omega_{2g-1} \cup \omega_{2g-1}]
        \end{pmatrix} \begin{pmatrix}
            [\omega_0 \cup \omega_0] & ... & [\omega_{2g-1} \cup \omega_0] \\ 
             \vdots & \ddots & \vdots \\
            [\omega_0 \cup \omega_{2g-1}] & ... & [\omega_{2g-1} \cup \omega_{2g-1}]
        \end{pmatrix} \\
        &= \begin{pmatrix}
            \sum\limits_{k=0}^{2g-1} [T_p\omega_0 \cup \omega_k][\omega_0 \cup \omega_k] & ... & \sum\limits_{k=0}^{2g-1} [T_p\omega_0 \cup \omega_k][\omega_{2g-1} \cup \omega_k] \\
            \vdots & \ddots & \vdots \\
            \sum\limits_{k=0}^{2g-1} [T_p\omega_{2g-1} \cup \omega_k][\omega_0 \cup \omega_k] & ... & \sum\limits_{k=0}^{2g-1} [T_p\omega_{2g-1} \cup \omega_k][\omega_{2g-1} \cup \omega_k]        
        \end{pmatrix}.
    \end{align*}
Therefore, for an entry of the matrix $(MC^{-1})^\intercal$, we have 
\begin{equation*}
    (MC^{-1})^\intercal _{i, j} = \sum\limits_{k=0}^{2g-1} [T_p\omega_j \cup \omega_k][\omega_i \cup \omega_k], \hspace{0.5cm}  \text{ where } [\omega_i \cup \omega_k] = \begin{cases}
        1 \hspace{0.8cm} \text{ if } k = i + g,\\
        -1 \hspace{0.5cm}\text{ if } i = k + g,\\
        0 \hspace{0.8cm}\text{ otherwise.}
    \end{cases}
\end{equation*}
On the other hand, if we have $T_p(\omega_i) = \sum\limits_{j = 0}^{2g-1} t_{i, j} \omega_j$, then $[T_p\omega_i \cup \omega_k] = \sum\limits_{j = 0}^{2g-1} t_{i, j} [\omega_j \cup \omega_k]$ and so 
\begin{align*}
    (MC^{-1})^\intercal _{i, j} &= \sum\limits_{k=0}^{2g-1} [\omega_i \cup \omega_k]\sum\limits_{\ell = 0}^{2g-1} t_{j, \ell} [\omega_\ell \cup \omega_k] \\
    &= \sum\limits_{k=0}^{g-1}[\omega_i \cup \omega_k][\omega_{k+g} \cup \omega_k] t_{j, k+g} + \sum\limits_{k=g}^{2g-1}[\omega_i \cup \omega_k][\omega_{k-g} \cup \omega_k] t_{j, k-g} \\
    &= -\sum\limits_{k=0}^{g-1}[\omega_i \cup \omega_k] t_{j, k+g} + \sum\limits_{k=g}^{2g-1}[\omega_i \cup \omega_k] t_{j, k-g} \\
    &= t_{j, i} \\
    &= (T_p)_{ij}.
\end{align*}
\end{proof}
\end{prop}

\par Therefore, by Proposition \ref{prop: Tp = (MC^-1)^T}, in order to compute the action of $T_p$ as a matrix, we need to compute the matrix $M$. This is equivalent to computing the action of $T_p$ on $q$-expansions of differentials and then computing cup products between the basis elements $\{\omega_0, ..., \omega_{2g-1}\}$ and their images under $T_p$, from the definition of the matrix $M$. The action of $T_p$ on $q$-expansions of modular forms, for $p \not| N$, is the following:
\begin{align*}
    T_p : \mathcal{S}_2^+(\Gamma_0(N)) &\longrightarrow \mathcal{S}_2^+(\Gamma_0(N)), \\ 
    \sum_{n = 0}^{\infty} a_nq^n &\mapsto \sum_{n = 0}^{\infty} a_{pn} q^n + p \sum_{n = 0}^{\infty} a_n q^{pn}.
\end{align*}
This induces a map on differentials of the form $\sum\limits_{n = -M}^{\infty} a_nq^n dq$, and so if we define $a_n \coloneqq 0$ for all $n < -M$, then we can write this induced map as
\begin{equation}
    \begin{aligned} \label{eqn: Tp induced map}
        T_p : \Omega_{X_0^+(N)} &\longrightarrow \Omega_{X_0^+(N)}, \\ 
    \sum_{n = -M}^{\infty} a_nq^n dq&\mapsto \sum_{n = -M}^{\infty} a_{pn} q^n dq+ p \sum_{n = -M}^{\infty} a_n q^{pn}dq.
    \end{aligned}
\end{equation}
By the definition of the map in (\ref{eqn: Tp induced map}), we note that $T_p$ sends a pole of order $k$ to a pole of order $pk$, therefore if a $q$-expansion $f$ has no poles then $T_p f$ will also have no poles. Therefore, the top left entry of the cup product matrix $M$ is zero. We choose $p$ so that $T_p$ does not lie in $\Z \subset \text{End}(J)$. Using the definition of the induced map of $T_p$ as in (\ref{eqn: Tp induced map}), we can calculate the image of $\text{H}^1_{\dR}(X/\Q)$ under the action of $T_p$, so that we have $T_p \omega_i$ for $i \in \{0, ..., 2g-1\}$ as $q$-expansions. We then use Serre's cup product formula to compute the cup product matrix $M$. Finally, we compute the matrix $T_p = (M C^{-1})^\intercal $. A benefit of our method is that we compute the entries of $T_p$ as exact rational numbers, rather than as a $p$-adic approximation for a given prime $p$, which allows us to work with multiple primes.

\par Next, we choose $Z$ to be a polynomial in $T_p$, and for $Z$ to be a nice correspondence it is enough for the corresponding matrix to have trace zero. We find that as in \cite{balakrishnan2023quadratic} and \cite{balakrishnan2021two}, we can take
\begin{equation*}
    Z = (\text{Tr}(T_p) I_{2g} - 2g T_p)C^{-1}
\end{equation*}
as our nice correspondence. We therefore can compute the matrix $\Lambda$ as in (\ref{def of lambda}). 

\par Once we compute the matrix $\Lambda$, we next need to calculate the matrix $G_\infty$, as defined in (\ref{matrix: Gx def}). We write $\psi_{\alpha \beta} \coloneqq G_\infty$, and so $\psi_{\alpha \beta}$ can be written as
\begin{equation*}
    \psi _{\alpha \beta} \coloneqq \begin{pmatrix}
        1 & 0 & 0 \\ \boldsymbol{a} & I_{2g} & 0 \\ c & \boldsymbol{b}^\intercal & 1 
    \end{pmatrix},
\end{equation*}
so that it satisfies (\ref{Gx and Lambda eqn}). Therefore we compute $\psi_{\alpha \beta}$ from this defining equation which is equivalent to requiring that
\begin{align*}
    d\boldsymbol{a} + \boldsymbol{\omega} &= 0, \\
    d\boldsymbol{b}^\intercal + \boldsymbol{\omega}^\intercal  Z &= 0, \\
    dc + \boldsymbol{b}^\intercal  \boldsymbol{\omega}  &= 0.
\end{align*}

\par We will now complete our computation of the Hodge filtration by computing $\boldsymbol{\beta}_{\text{Fil}}, \gamma_{\text{Fil}}$ as discussed in earlier sections. To do this we will follow Algorithm 5.20 in \cite{balakrishnan2023computation}. At this point, we have calculated the matrix $\psi_{\alpha \beta}$ so all that is left is to compute $b_{\text{Fil}} = (b_g, ..., b_{2g-1}) \in \Q_p^g$ and the function $\gamma_{\text{Fil}}$ characterised by $\gamma_{\text{Fil}}(b) = 0$ such that 
\begin{equation*}
    g_\infty + \gamma_{\text{Fil}} - b_{\text{Fil}}^\intercal N^\intercal \boldsymbol{\Omega}^\intercal - \boldsymbol{\Omega}^\intercal ZNN^\intercal \boldsymbol{\Omega} \in K[[t]],
\end{equation*}
where we write $N = (0_g, I_g)^\intercal \in M_{2g \times g}(\Q)$. The existence and uniqueness of $\boldsymbol{\beta}_{\text{Fil}}, \gamma_{\text{Fil}}$ follows from the following. We first recall Lemma 4.7 from \cite{balakrishnan2019explicit}, with our notation:

\begin{lemma}

Given any $g \in \Q((t))$, there exists a unique vector of constants  $\beta \in \Q_p^{2g}$ and a function $\gamma \in \text{H}^0(Y, \mathcal{O})$, unique modulo constants, such that, 

\begin{equation*}
g - \gamma_{\text{Fil}} - b_{\text{Fil}}^\intercal N^\intercal \boldsymbol{\Omega}^\intercal - \boldsymbol{\Omega}^\intercal Z N N^\intercal \boldsymbol{\Omega}^\intercal \in \Q[[t]].
\end{equation*}

\end{lemma}
 Set $\boldsymbol{\beta}_{\text{Fil}} = (0, ..., b_g, ..., b_{2g-1})^\intercal$. Rewriting this expression in the notation of $\psi_{\alpha \beta}$ (since $G_\infty = \psi_{\alpha \beta}$), we will solve for  $b_{\text{Fil}} = (b_g, ..., b_{2g-1}) \in \Q_p^g$ and the $q$-expansion $\gamma_{\text{Fil}}$ such that
\begin{equation} \label{eqn: beta fil and gamma fil expression}
    c + \gamma_{\text{Fil}} - b_{\text{Fil}} N^\intercal \boldsymbol{a} - \boldsymbol{b}^\intercal NN^\intercal \boldsymbol{a} \in \Q[[q]].
\end{equation}

\par We claim that we are able to construct the function $\gamma_{\text{Fil}}$ as a linear combination of elements of the set $\{ \frac{s_{d-(j_{\dR} - 1)}}{s_d}, ..., \frac{s_{d-1}}{s_d} \}$. We first bound the order of the maximum pole of $\{\int \omega_i\}$ for $i \in \{0, ..., 2g-1\}$ in terms of $j_{\dR} \in \mathbb{N}$. Let $-M$ be the maximum pole order of $\{\int \omega_i\}$ for $i \in \{0, ..., 2g-1\}$, with $M > 0$, and we recall that $\{\int \omega_i\}$ for $i \in \{0, ..., 2g-1\}$ corresponds to the rows of $\boldsymbol{a}$. Then, since $\{\omega_0, ..., \omega_{g-1}\}$ has increasing valuation and $val(\omega_0) = 0$, we have
\begin{equation*}
    -M = val\left(\int f_{\dR} \cdot \omega_0\right) = val(s_{d - j_{\dR}}) - val(s_d) + 1.
\end{equation*}
Since $\{s_1, ..., s_d\}$ is a basis with strictly increasing valuation and $val(s_1) = 1$, we have $val(s_{d - j_{\dR}}) \geq d - j_{\dR}$. Moreover, $val(s_d) = \frac{N+1}{2}$ and we have $\frac{N+1}{2} = d + g$ by the equality in (\ref{eqn: (N+1)/2 = d + g}). As a result, we have 
\begin{equation*}
    M \leq g + j_{\dR} - 1.
\end{equation*}
Therefore, the expression $c - \boldsymbol{b}^\intercal NN^\intercal \boldsymbol{a}$ has a pole of order at most $M$. Therefore, the sets $\{ \frac{s_{d-(j_{\dR} - 1)}}{s_d}, $ $ ..., \frac{s_{d-1}}{s_d} \}$ and $\{ \int \omega_g, ..., \int \omega_{2g-1}\}$ give $j_{\dR} - 1 + g$ elements with principal parts spanning the space up to the pole order of the expression $c - \boldsymbol{b}^\intercal NN^\intercal \boldsymbol{a}$, which is $q^{-M}\Q((q))/\Q[[q]]$.

\par All that is left is to show that the sets $\{ \frac{s_{d-(j_{\dR} - 1)}}{s_d}, ..., \frac{s_{d-1}}{s_d} \}$ and $\{ \int \omega_g, ..., \int \omega_{2g-1}\}$ are linearly independent from each other. Firstly, the sets are individually linearly independent because the elements of $\{ \frac{s_{d-(j_{\dR} - 1)}}{s_d}, ..., \frac{s_{d-1}}{s_d} \}$ all have different valuations by the definition of the basis $\{s_1, ..., s_d\}$. Similarly, the set $\{ \int \omega_g, ..., \int \omega_{2g-1}\}$ is linearly independent because it is built from the basis $\{\omega_0, ..., \omega_{g-1}\}$ as outlined in \S \ref{subsec: computing H1dR}. Now, suppose there is a linear dependence from these sets, so suppose there exist constants $\lambda_i$ for $i = g, ..., 2g-1$ and $\mu_i$ for $i = 1, ..., j_{\dR} - 1$ such that 
\begin{equation} \label{eqn: mu and lambda coefficients eqn}
    \sum\limits_{i=1}^{j_{\dR} - 1} \mu_i \frac{s_{d-i}}{s_d} + \sum\limits_{i=g}^{2g-1} \lambda_i \int \omega_i = 0
\end{equation}
and not all the $\lambda_i, \mu_i$ are zero. Taking the derivative of (\ref{eqn: mu and lambda coefficients eqn}) and writing 
\begin{equation*}
    f \coloneqq \sum\limits_{i=1}^{j_{\dR} - 1} \mu_i \frac{s_{d-i}}{s_d} \text{ and } \omega \coloneqq \sum\limits_{i=g}^{2g-1} \lambda_i \omega_i,
\end{equation*}
gives that the expression $df + \omega = 0$ has no poles on the affine open $Y = X \setminus \{\infty\}$, since $s_d$ and $\omega_g, ..., \omega_{2g-1}$ only have poles at $\infty$. Therefore, the function $f$ is constant and $\omega \in \text{H}^0(X, \Omega^1)$ by \cite{balakrishnan2021quadratic} Lemma 6.3. It follows that we have $\mu_i = 0$ for all $i = 1, ..., j_{\dR} - 1$, and $\lambda_i = 0$ for all $i = g, ..., 2g-1$, and so the sets are linearly independent. Therefore, we can construct the function $\gamma_{\text{Fil}}$ from the set $\{ \frac{s_{d-(j_{\dR} - 1)}}{s_d}, ..., \frac{s_{d-1}}{s_d} \}$.

\par Once we have constructed $\gamma_{\text{Fil}}$ so that (\ref{eqn: beta fil and gamma fil expression}) holds, $\gamma_{\text{Fil}}$ is then determined uniquely by the condition that $\gamma_{\text{Fil}}(b) = 0$ for $b \in X(\Q)$. Since we construct $\gamma_{\text{Fil}}$ from the set $\{ \frac{s_{d-(j_{\dR} - 1)}}{s_d}, ..., \frac{s_{d-1}}{s_d} \}$, the elements of which have poles only at $\infty$, we can take $b$ to be the upper half plane representative of a rational CM point of $X$.

\par In summary, in this section we have developed the following algorithm.

\begin{algorithm}[Model-free Computation of the Hodge filtration for $X_0^+(N)$ for $N$ prime] \label{alg: model-free HF}
\hfill \break 
    Input: 
    \begin{itemize}
        \item A prime $N$ such that $X_0^+(N)/\Q$ is a modular curve such that $g > 1$.
        \item A prime $p$ such that $T_p$ does not lie in $\Z \subset \text{End}(J)$.
    \end{itemize}

    Output: the matrix $\lambda^{\text{Fil}}$, defined in terms of weakly holomorphic modular forms.

    Steps:
    \begin{enumerate}
        \item Compute a basis of $q$-expansions of $\mathcal{S}_2^+(\Gamma_0(N))$
        \item Compute a basis of $q$-expansions of $\mathcal{S}_{12}^+(\Gamma_0(N))$, denoted by $\{s_1, ..., s_d\}$ such that $val(s_i) < val(s_{i+1})$ for all $i \in \{1, ..., d-1\}$ where $d = \dim \mathcal{S}_2^+(\Gamma_0(N))$.
        \item Compute $j_{\dR} \in \{1, ..., d-1\}$ to be the smallest integer such that (\ref{eqn: condition on j}) holds.
        \item Compute the matrices $T_p, Z$.
        \item Compute a vector $\boldsymbol{\Omega}$ that satisfies $d\boldsymbol{\Omega} = -\boldsymbol{\omega}$.
        \item Solve the system of equations, (\ref{Gx and Lambda eqn}), to determine the matrix $\psi_{\alpha, \beta}$.
        \item Compute the vector of constants $\boldsymbol{\beta}_{\text{Fil}} \in \Q_p^{2g}$, and $\gamma_{\text{Fil}} \in span\{\frac{s_{d-(j_{\dR} - 1)}}{s_d}, ..., \frac{s_{d-1}}{s_d}\}$ such that (\ref{eqn: beta fil and gamma fil expression} holds), and such that $\gamma_{\text{Fil}}(b) = 0$.
    \end{enumerate}

\end{algorithm}

\par Once we have computed the Hodge filtration, we will investigate congruences of coefficients of modular forms modulo $N$, in particular, the potential congruences between iterated integrals in the plus eigenspace for the Atkin-Lehner involution and single integrals in the minus eigenspace for the Atkin-Lehner involution. Let $g_0(N)$ be the genus of $X_0(N)$, and let $\{\tilde{s}_0, ..., \tilde{s}_{g_0(N) - 1}\}$ be a basis of $q$-expansions for $\mathcal{S}_2(\Gamma_0(N))$, so that $\{\tilde{\omega}_0, ..., \tilde{\omega}_{g_0(N) - 1}\}$ is a basis of $q$-expansions for $\text{H}^0(X_0(N), \Omega^1)$ where we let $\tilde{\omega}_i \coloneqq \frac{\tilde{s}_i}{q}dq$ for all $i \in \{0, ..., g_0(N) - 1\}$. Define the following sets
\begin{align*}
    A_1 &\coloneqq \left\{\left(\int \omega_i \right)\cdot \omega_j: \text{ for } i, j, \in \{0, ..., 2g-1\}\right\}, \\
    A_2 & \coloneqq  \left\{ \tilde{\omega_i}: \text{ for } i \in \{0, ..., g_0(N)-1\} \right\}, \\
    A_3 & \coloneqq  \left\{ \omega_i: \text{ for } i \in \{0, ..., g-1\} \right\}.
\end{align*}

\par Choose the positive integer $n_N$, and for two sets $A, A'$, construct a matrix so that each row is the first $n_{N}$ coefficients of an element of $A \cup A'$, with each coefficient taken modulo $N$. Denote such a matrix by $M_N(A \cup A')$. We will compute the matrices $M_N(A_1 \cup A_2), M_N(A_1 \cup A_3)$ and then compute their coranks. We note that these matrices have $4g^2+g_0(N), 4g^2 + g$ rows respectively. We will then compute the value 
\begin{equation} \label{eqn: corank quantity}
corank(M_N(A_1 \cup A_2)) - corank(M_N(A_1 \cup A_3)).
\end{equation}
This gives the corank of the matrix of coefficients of elements in $A_1$, which lie in the plus eigenspace for the Atkin-Lehner involution, and the elements in $A_2 \setminus A_3$, which lie in the minus eigenspace for the Atkin-Lehner involution. Therefore, if the value in (\ref{eqn: corank quantity}) is nonzero, this suggests a relation between the elements of $A_1$ and elements of the form $\{ \frac{s}{q} dq: s \in \mathcal{S}_2^{-}(\Gamma_0(N)) \}$. Moreover, this means that the integrals of these elements also have a relation, thus showing there is a relation between the iterated integrals on the plus eigenspace, and the single integrals on the minus eigenspace. 


\section{Examples}\label{sec: example}

\par In this section, we apply our algorithm to compute Hodge filtrations on the modular curves $X_0^+(67)$, $ X_0^+(193)$, of genus 2 and 7, respectively. We also compute congruences of modular forms for the genus 3 curve $X_0^+(97)$.

\subsection{Genus 2 example}

\par We choose the curve $X_0^+(67)$, since $N=67$ is the first prime such that the curve $X_0^+(N)$ or prime level has genus 2. The 10 rational points were computed by Balakrishnan-Best-Bianchi-
Lawrence-M\"uller-Triantafillou-Vonk \cite{balakrishnan2021two} using the quadratic 
Chabauty algorithm and so requiring an explicit plane model of the curve.

\par Following \S \ref{subsec: computing H1dR}, the first step of implementing quadratic Chabauty is to compute a symplectic basis of $\text{H}^1_{\dR}(X/\Q)$. By computing a basis of $\mathcal{S}_2^+(\Gamma_0(67))$, denoted by $\{f_{\omega_0}, f_{\omega_1} \}$, we obtain the following first two basis elements:
\begin{align*}
    f_{\omega_0} &= q - 3q^3 - 3q^4 - 3q^5 + O(q^6) &\mapsto \hspace{0.6cm} \omega_0 &= (1 - 3q^2 - 3q^3 -3q^4 + O(q^5)) dq, \\
    f_{\omega_1} &= q^2 - q^3 - 3q^4 + O(q^6) &\mapsto  \hspace{0.6cm} \omega_1 &= (q - q^2 - 3q^3 + O(q^5)) dq.
\end{align*}    
We have $d = \dim \mathcal{S}_{12}^+(\Gamma_0(N)) = 32$, we compute a basis of increasing valuation $\{s_1, ..., s_{32}\}$, and check that $s_{32} = a_{\Delta} \sqrt{\Delta \cdot w_N(\Delta)}$. We find that $j=1$ satisfies (\ref{eqn: condition on j}), since we compute the following
\begin{align*}
    s_{31} &= q^{31} - 11q^{32} + 45q^{33} - 856q^{35} + 3862q^{36} - 5982q^{37} -
    6709q^{38} + 35355q^{39} + O(q^{40}), \\
    s_{32} &= q^{34} - 12q^{35} + 54q^{36} - 88q^{37} - 99q^{38} + 540q^{39} + O(q^{40}),
\end{align*}
and $\max\limits_{i=0, 1} (val(f_{\omega_i})) = 2$. Therefore, we take their quotient to give our function with poles only at $\infty$ as 
\begin{equation*}
    f_{\dR} \coloneqq \frac{s_{31}}{s_{32}} = q^{-3} + q^{-2} + 3q^{-1} + 70 + 9q + O(q^2).
\end{equation*}
Write $\eta_{i+2} \coloneqq f_{\dR}\cdot \omega_i$ for $i=0, 1$. We compute the cup product of our four $q$-expansions $\{\omega_0, \omega_1, \eta_2, \eta_3\}$ and obtain the following
\begin{equation*}
    \begin{pmatrix}
        0 & 0 & -1 & -1 \\ 0 & 0 & - \frac{1}{2} & 0 \\ 1 & \frac{1}{2} & 0 & -33 \\ 1 & 0 & 33 & 0 
    \end{pmatrix}.
\end{equation*}
Since this has nonzero determinant, we have verified that we have indeed found a basis of $\text{H}^1_{\dR}(X/\Q)$ consisting of $q$-expansions. We now want to construct a symplectic basis $\{\omega_0, \omega_1, \omega_2, \omega_3\}$, so that its cup product matrix is equal to $C$. We will do this using a change of basis to construct the two new final basis elements $\{ \omega_2, \omega_3 \}$ as
\begin{align*}
    &\omega_2 \coloneqq -66\omega_1 - \eta_3 = (-q^{-2} + 1 - 130q + 136q^2 + 401q^3 + 10q^4 + 13q^5 + O(q^6))dq, \\
    &\omega_3 \coloneqq -2\eta_2 + 2\eta_3 = (-2q^{-3} - 130 + 140q + 276q^2 + 30q^3 + 446q^4 -
        40q^5 + O(q^6))dq.
\end{align*} 
This results in $\{\omega_0, \omega_1, \omega_2, \omega_3\}$ having $C$ as in (\ref{def: symp cup product matrix}) as its cup product matrix and therefore is a symplectic basis of $\text{H}^1_{\dR}(X/\Q)$. Next, following the method to calculate $T_p$, with $p=3$, we find that
\begin{equation*}
    M = \begin{pmatrix}
        0  &  0  &  3  & -1 \\ 
  0  &  0  &  1  &  0 \\ 
  -3 &  -1  &  0 & 263 \\ 
   1  &  0 & -263  &  0
    \end{pmatrix}, \hspace{1cm} T_3 = \begin{pmatrix}
    3  &  1  &  0 & -263 \\
  -1  &  0 & 263  &  0 \\ 
   0  &  0   & 3  & -1 \\ 
   0  &  0  &  1  &  0
    \end{pmatrix}.
\end{equation*}
Therefore, we obtain the following for $Z$
\begin{equation*} 
    Z = (\text{Tr}(T_3) I_{2g} - 2g T_3)C^{-1} = \begin{pmatrix}
        0 & 1052  &  6   & 4 \\
 -1052  &   0  &   -4  &   -6 \\
    -6  &  4  &   0   &  0 \\
    -4  &  6  &   0   &  0
    \end{pmatrix}.
\end{equation*}

Since we have $\eta = 0$, we obtain the following for the matrix $\Lambda$, as defined in (\ref{def of lambda})

\noindent\resizebox{\columnwidth}{!}{
$
    \begin{pmatrix}
        0 & 0 & 0 & 0 & 0 & 0 \\[6pt] 
        1 - 3q^2 +O(q^3) & 0 & 0 & 0 & 0 & 0 \\[6pt] 
        q - q^2 + O(q^3) & 0 & 0 & 0 & 0 & 0 \\[6pt] 
        -q^{-2} + 1 + O(q) & 0 & 0 & 0 & 0 & 0 \\[6pt] 
        -2q^{-3} -  130 + O(q) & 0 & 0 & 0 & 0 & 0 \\[6pt] 
        0 & 8q^{-3} + 6q^{-2} + O(1) & -12q^{-3} - 4q^{-2} + O(1) 
 & 6 - 4q + O(q^2) & 4 - 6q + O(q^2) & 0 \\[6pt] 
    \end{pmatrix} \textit{dq}.
$
}

Solving the relation (\ref{Gx and Lambda eqn}) (and for now setting all constants of integration to 0) we obtain the following for the matrix $\psi _{\alpha \beta}$

\noindent\resizebox{\columnwidth}{!}{%
$
    \begin{pmatrix}
        1 & 0 & 0 & 0 & 0 & 0\\[6pt] 
        -q + q^3 + O(q^4) & 1 & 0 & 0 & 0 & 0 \\[6pt] 
    - \frac{1}{2}q{^2} + \frac{1}{3}q^3 + O(q^4) & 0 & 1 & 0 & 0 & 0 \\[6pt] 
    -q^{-1} - q  + O(q^2) & 0 & 0 & 1 & 0 & 0 \\[6pt] 
    -q^{-2} + 130q  + O(q^2) & 0 & 0 & 0 & 1 & 0 \\[6pt] 
    12q^{-1} + 16q  + O(1)
 & 4q^{-2} + 6q^{-1} + O(1) &
    - 6q^{-2} - 4q^{-1} + O(1) &
    - 6q + 2q^2 + O(q^3) &
    - 4q + 3q^2 + O(q^3) & 1
    \end{pmatrix}.
$
}

\par We can now complete our computation of the Hodge filtration by computing $\boldsymbol{\beta}_{\text{Fil}}, \gamma_{\text{Fil}}$ as discussed in §\ref{subsec: calculating hodge}. Let us write $\boldsymbol{a} = (A_1, A_2)^\intercal, \boldsymbol{b} = (B_1, B_2)^\intercal$, so that they are a tuple of 2-dimensional vectors. By the definition of the matrix $N$, we would like the following expression to have no poles:
\begin{equation*}
    c + \gamma_{\text{Fil}} - b_{\text{Fil}} A_2 - B_2^ \intercal A_2.
\end{equation*}
We have already calculated $c, A_2, B_2$ and we find that $c - B_2^ \intercal A_2 = 8q^{-1} - 3 + O(q)$. Therefore, using the fact that $A_2 = \begin{pmatrix}
    -q^{-1} + O(q) & -q^{-2} + O(q)
\end{pmatrix}$, we can cancel this pole of by setting $b_{\text{Fil}} = \begin{pmatrix}
    -8 & 0
\end{pmatrix}$. Additionally, since $j_{\dR} = 1$, we have $\{ \frac{s_{d - (j_{\dR} - 1)}}{s_d}, ..., \frac{s_{d-1}}{s_d} \} = \emptyset$, and so $\gamma_{\text{Fil}} = 0$. Hence we have calculated the Hodge filtration ending with
\begin{equation*}
    \boldsymbol{\beta}_{\text{Fil}} = \begin{pmatrix}
        0 & 0 & -8 & 0 
    \end{pmatrix}^\intercal, \hspace{0.5cm} \gamma_{\text{Fil}} = 0.
\end{equation*}

\par  We finish this example by investigating congruences of modular forms. For $N=67$, take $n_{67} = 60$, and we find that the quantity (\ref{eqn: corank quantity}) is 
\begin{align*}
    corank(M_{67}(A_1 \cup A_2)) - corank(M_{67}(A_1 \cup A_3)) &= 3 - 2 = g - 1.
\end{align*}
The fact that the difference of coranks of these matrices is positive suggests that there is a relationship between the iterated integrals on the plus eigenspace, and the single integrals on the minus eigenspace. Moreover, by Remark \ref{remark: rho = g} we notice that the difference is equal to $\rho(J) - 1$ which is equal the rank of $\Ker(NS(J) \rightarrow NS(X))$ from which we construct $Z$.

\subsection{Genus 3 example}

\par The prime $N = 97$ is the smallest integer such that $X_0^+(N)$ is a nice curve of genus 3. We will investigate congruences of integrals for this curve. We being with computing a basis of $\text{H}^1_{\text{dR}}(X/\Q)$, following \S \ref{subsec: computing H1dR}. By computing a basis of $\mathcal{S}_2^+(\Gamma_0(97))$, denoted by $\{f_{\omega_0}, f_{\omega_1}, f_{\omega_2}\}$, we obtain the following the first three basis elements:
\begin{align*}
    f_{\omega_0} &= q - 4q^4 - 5q^5 - 3q^6 - q^7 + O(q^8)  &\mapsto \hspace{0.6cm} \omega_0 &= (1 - 4q^3 - 5q^4 - 3q^5 - q^6 + O(q^7))dq,  \\
    f_{\omega_1} &= q^2 - 3q^4 - q^5 - 2q^6 + O(q^8)  &\mapsto \hspace{0.6cm} \omega_1 &= (q - 3q^3 - q^4 - 2q^5 + O(q^7))dq,  \\
    f_{\omega_2} &= q^3 - q^4 - 2q^5 - q^6 + q^7 + O(q^8) &\mapsto \hspace{0.6cm} \omega_2 &= (q^2 - q^3 - 2q^4 - q^5 + q^6 + O(q^7))dq.
\end{align*}
We have $d = \dim \mathcal{S}_{12}^+(\Gamma_0(97)) = 46$, we compute a basis of increasing valuation $\{s_1, ..., s_{46}\}$, and check that $s_{46} = a_\Delta \sqrt{\Delta \cdot w_N(\Delta)}$. We find that $j=1$ satisfies (\ref{eqn: condition on j}), and so we take their quotient to give our function with poles only at $\infty$ as 
\begin{equation*}
    f_{\text{dR}} \coloneqq \frac{s_{45}}{s_{46}} = q^{-4} + q^{-3} + 3q^{-2} + 4q^{-1} + 73 + O(q). 
\end{equation*}
Write $\eta_{i + 3} \coloneqq f_{\text{dR}} \cdot \omega_i$ for $i = 0 , 1, 2$. We compute the cup product matrix of our 6 $q$-expansions $\{\omega_0, \omega_1, \omega_2,$ $ \eta_3, \eta_4, \eta_5\}$ and find that it has nonzero determinant we have indeed found a basis of $\text{H}^1_{\text{dR}}(X/\Q)$. Taking a change of basis we obtain the following final 3 basis elements to give a symplectic basis denoted by $\{\omega_0, ..., \omega_5\}$
\begin{align*}
    \omega_3 &\coloneqq -2\omega_1 - 129\omega_2 - \eta_5,   \\
    \omega_4 &\coloneqq  60\omega_2 - 2\eta_4 + 2\eta_5,  \\
    \omega_5 &\coloneqq -3\eta_3 + 3\eta_4 + 6\eta_5.   
\end{align*}
Now that we have computed the symplectic basis, we choose $n_{97} = 90$, and compute the difference of coranks as
\begin{equation*}
    corank(M_{97}(A_1 \cup A_2)) - corank(M_{97}(A_1 \cup A_3)) = 11 - 9 = g - 1.
\end{equation*}
The fact that the difference of coranks of these matrices is positive suggests that there is a relationship between the iterated integrals on the plus eigenspace, and the single integrals on the minus eigenspace. Moreover, by Remark \ref{remark: rho = g} we notice that the difference is equal to $\rho(J) - 1$ which is equal the rank of $\Ker(NS(J) \rightarrow NS(X))$ from which we construct $Z$.

\subsection{Genus 7 example}
\par We begin by computing the smallest prime $N$ such that the genus of $X_0^+(N)$ is 7. Denote the genus of $X_0^+(N)$ by $g_0^+(N)$. Following \S 4 of \cite{advzaga2021quadratic}, we compute that for a prime $N$ such that $g_0^+(N) \leq 7$ we must have $N \leq 13400$. Computing the genus, $g_0^+(N)$, for all $N \leq 13400$, we have that the genus of $X$ is 7 if and only if 
\begin{equation*}
    N \in \{193, 229, 233, 241, 257, 281 \}.
\end{equation*}
Therefore, we will carry out our method on $X_0^+(193)$. We begin with computing a basis of $\text{H}^1_{\dR}(X/\Q)$, following §\ref{subsec: computing H1dR}. By computing a basis of $\mathcal{S}_2^+(\Gamma_0(193))$, denoted by $\{f_{\omega_0}, ...,f_{\omega_6} \}$, we obtain the following first seven basis elements:
\begin{align*}
    f_{\omega_0} &= q - 6q^6 - 7q^8 - 4q^{10} + O(q^{11}) &\mapsto \hspace{0.6cm} \omega_0 &= (1 - 6q^5 - 7q^7 - 4q^9 + O(q^{10})) dq, \\
    f_{\omega_1} &= q^2 - 5q^6 - 3q^8 - q^{10} + O(q^{11}) &\mapsto \hspace{0.6cm} \omega_1 &= (q - 5q^5 - 3q^7 - q^9 + O(q^{10})) dq, \\
    f_{\omega_2} &= q^3 - 3q^6 - 2q^8 - q^{10} + O(q^{11}) &\mapsto \hspace{0.6cm} \omega_2 &= (q^2 - 3q^5 - 2q^7 - q^9 + O(q^{10})) dq, \\
    f_{\omega_3} &= q^4 - 2q^6 - 2q^8 + q^{10} + O(q^{11}) &\mapsto \hspace{0.6cm} \omega_3 &= (q^3 - 2q^5 - 2q^7 + q^9 + O(q^{10})) dq, \\
    f_{\omega_4} &= q^5 - q^6 - q^8 - q^{10} + O(q^{11}) &\mapsto \hspace{0.6cm} \omega_4 &= (q^4 - q^5 - q^7 - q^9 + O(q^{10})) dq, \\
    f_{\omega_5} &= q^7 - q^8 - q^{10} - 6q^{11} + O(q^{12}) &\mapsto \hspace{0.6cm} \omega_5 &= (q^6 - q^7 - q^9 - 6q^{10} + O(q^{11})) dq, \\
    f_{\omega_6} &= q^9 - q^{10} - 3q^{11} + q^{12} + O(q^{13}) &\mapsto \hspace{0.6cm} \omega_6 &= (q^8 - q^9 - 3q^{10} + q^{11}+ O(q^{12})) dq.
\end{align*}  
We have $d = \dim \mathcal{S}_{12}^+(\Gamma_0(N)) = 90$, we compute a basis of increasing valuation $\{s_1, ..., s_{90}\}$, and check that $s_{90} = a_\Delta \sqrt{\Delta \cdot w_N(\Delta)}$. We find that $j = 3$ is the smallest $j_{\dR} \in \mathbb{N}$ that satisfies (\ref{eqn: condition on j}), since we compute the following
\begin{align*}
    s_{87} &= q^{87} - 11q^{88} + 429q^{90} - 19711q^{92} + 94667q^{93} - 128929q^{94} -
    O(q^{95}), \\
    s_{90} &= q^{97} - 12q^{98} + 54q^{99} - 88q^{100} - 99q^{101} + 540q^{102} -
    418q^{103} - 648q^{104} + O(q^{105}).
\end{align*}
and $\max\limits_{i=0, ..., 6} (val(f_{\omega_i})) = 9$. Therefore, we take their quotient to give our function with poles only at $\infty$ as 
\begin{multline*}
    f_{\dR} \coloneqq \frac{s_{87}}{s_{90}} = q^{-10} + q^{-9} - 42q^{-8} - 41q^{-7} + 1963q^{-6} + 1922q^{-5} + 3841q^{-4} \\ + 5806q^{-3} +
    9690q^{-2} + 11488q^{-1} + 284925 + O(q).
\end{multline*}
Write $\eta_{i+7} \coloneqq f_{\dR}\cdot \omega_i$ for $i=0, ..., 6$. We compute the cup product of our fourteen $q$-expansions $\{\omega_0, ..., \omega_6,$ $ \eta_7, ..., \eta_{13} \}$ and obtain the following matrix

\noindent\resizebox{\columnwidth}{!}{
$
\begin{pmatrix}
        0  & 0 & 0 & 0 & 0 & 0 & 0 & \frac{-79753}{8} & \frac{-48393}{8} & -3961 & -2001 & -2003 &  43 & -1 \\[6pt]
0  & 0 & 0 & 0 & 0 & 0 & 0 & \frac{-73423}{24}  & \frac{-16469}{8}  & \frac{-6133}{6}  & \frac{-3065}{3}  & 0 &  0  & 0 \\[6pt]
0  & 0 & 0 & 0 & 0 & 0 & 0 & \frac{-15583}{12}   & \frac{-2639}{4}   & \frac{-3917}{6} &  \frac{89}{6}  & \frac{43}{3}   & \frac{-1}{3}  &  0 \\[6pt]
0  & 0 & 0 & 0 & 0 & 0 & 0 & \frac{-5909}{12} & \frac{-2013}{4}  & \frac{34}{3}  &  \frac{34}{3}  & 0 & 0 & 0 \\[6pt]
0  & 0 & 0 & 0 & 0 & 0 & 0 & \frac{-47917}{120}  & \frac{53}{40}  & \frac{257}{30} & \frac{-1}{30}  & \frac{-1}{5}  & 0 & 0\\[6pt]
0  & 0 & 0 & 0 & 0 & 0 & 0 & \frac{49}{8}  &  \frac{-1}{56}  &  \frac{-1}{7} & 0 & 0 & 0 & 0 \\[6pt]
0  & 0 & 0 & 0 & 0 & 0 & 0 & \frac{-1}{9}  & 0 & 0 & 0 & 0 & 0 & 0  \\[6pt]
\frac{79753}{8} & \frac{73423}{24} & \frac{15583}{12} & \frac{5909}{12} & \frac{47917}{120} & \frac{-49}{8} & \frac{1}{9} & 0 & \frac{-70103211175}{84} & \frac{-15086097454}{21}  &   \frac{-24213252503}{60}  & \frac{-3856217914}{9}  &  \frac{1301390974}{45}  & \frac{3928216}{9}  \\[6pt]
\frac{48393}{8}  & \frac{16469}{8}  & \frac{2639}{4}  &  \frac{2013}{4}  & \frac{-53}{40}  & \frac{1}{56}  & 0 & \frac{70103211175}{84}  & 0 & \frac{-583720821}{7}  & \frac{-7976234419}{60} &  \frac{303227023}{42}  &  \frac{2879807041}{210}  & \frac{16160411}{14}  \\[6pt]
3961 & \frac{6133}{6}  & \frac{3917}{6}  & \frac{-34}{3}  & \frac{-257}{30} &  \frac{1}{7}  & 0 &  \frac{15086097454}{21}  & \frac{583720821}{7}  & 0 &  \frac{17263481}{30}  & \frac{128613395}{42}  &  \frac{2008835527}{210}  & \frac{8685898}{7}  \\[6pt]
2001 &  \frac{3065}{3}  & \frac{-89}{6}  &  \frac{-34}{3}   & \frac{1}{30}   & 0 & 0 & \frac{24213252503}{60}  & \frac{7976234419}{60}  & \frac{-17263481}{30}  &  0  & \frac{-2743136}{15}  & \frac{121447697}{30}  &  \frac{-764545}{6}  \\[6pt] 
2003 &  0 &  \frac{-43}{3}  & 0  & \frac{1}{5}   & 0 &  0 &  \frac{3856217914}{9}  &  \frac{-303227023}{42}  &  \frac{-128613395}{42}  &  \frac{2743136}{15}  &  0  & \frac{60737216}{15} &  \frac{-376507}{3}  \\[6pt]
-43 &  0  & \frac{1}{3}  &  0  & 0  & 0 &  0  & \frac{-1301390974}{45}  &  \frac{-2879807041}{210}   & \frac{-2008835527}{210}   & \frac{-121447697}{30}  &  \frac{-60737216}{15}  &  0 &  \frac{2077}{3}  \\[6pt]
1 & 0 & 0 & 0 & 0 & 0 & 0 & \frac{-3928216}{9}  &  \frac{-16160411}{14}   & \frac{-8685898}{7}  &  \frac{764545}{6}   & \frac{376507}{3}   & \frac{-2077}{3}   & 0
    \end{pmatrix}.
$
}

\par Since this cup product matrix has nonzero determinant, we have verified that we have indeed found a basis of $\text{H}^1_{\dR}(X/\Q)$ consisting of $q$-expansions. Via a change of basis, we obtain a symplectic basis of $\text{H}^1_{\dR}(X/\Q)$, denoted by $\{\omega_0, ..., \omega_{2g-1}\}$ as the following: 

\begin{align*}
    \omega_{7} &\coloneqq - \frac{354960083}{3 \cdot 23 \cdot 101 \cdot 617} \omega_1 + 2077 \omega_2 - \frac{2792316815}{2 \cdot 3 \cdot  23 \cdot 101 \cdot 617} \omega_3 -478660 \omega_4 \\ &\hspace{2cm} - \frac{124457164816556}{3 \cdot 23 \cdot 101 \cdot 617} \omega_5 + \frac{619199026843193}{2^2 \cdot  23 \cdot 101 \cdot 617} \omega_6 - \eta_{13}, \\
    \omega_{8} &\coloneqq - \frac{46557376937}{3^2 \cdot  23 \cdot 101 \cdot 617} \omega_2 - \frac{12832634963}{2 \cdot 3^2 \cdot 23 \cdot 101 \cdot 617} \omega_3 + \frac{85816262590}{3^2 \cdot  23 \cdot 101 \cdot 617} \omega_4 \\ &\hspace{2cm} - \frac{8564991957794}{3^2 \cdot  23 \cdot 101 \cdot 617} \omega_5 - \frac{1731576379443407}{2^2 \cdot 3^2 \cdot  23 \cdot 101 \cdot 617} \omega_6 - \frac{272}{3^2 \cdot  23 \cdot 101 \cdot 617}\eta_8 \\ &\hspace{2cm} +
\frac{34}{3^2 \cdot  23 \cdot 101 \cdot 617}\eta_9 - \frac{12112}{3^2 \cdot  23 \cdot 101 \cdot 617}\eta_{10} + \frac{1673}{3^2 \cdot  23 \cdot 101 \cdot 617}\eta_{11} \\ &\hspace{4cm} + \frac{1574}{3 \cdot 23 \cdot 101 \cdot 617}\eta_{12} + \frac{22598827}{3^2 \cdot  23 \cdot 101 \cdot 617}\eta_{13},  \\
    \omega_{9} &\coloneqq  \frac{4375157750404}{3^2 \cdot  23 \cdot 101 \cdot 617}\omega_3 + 589366\omega_4 + \frac{18257845587344468}{3^2 \cdot  23 \cdot 101 \cdot 617}\omega_5 +
\frac{108609257490452313}{23 \cdot 101 \cdot 617}\omega_6 \\ &\hspace{2cm} -3\eta_{12} -129\eta_{13}, \\
    \omega_{10} &\coloneqq  \frac{35663178459305}{2^2 \cdot 3^2 \cdot  23 \cdot 101 \cdot 617} \omega_4 + \frac{564892475732921}{2^2 \cdot 3 \cdot 23 \cdot 101 \cdot 617}\omega_5 +
\frac{182812738838122147}{2^3 \cdot 3^2 \cdot 23 \cdot 101 \cdot 617}\omega_6 \\ &\hspace{2cm} - \frac{24520}{3^2 \cdot  23 \cdot 101 \cdot 617}\eta_8 + \frac{3065}{3^2 \cdot  23 \cdot 101 \cdot 617}\eta_9 +
\frac{92681}{2 \cdot 3^2 \cdot 23 \cdot 101 \cdot 617}\eta_{10} \\ &\hspace{2cm} - \frac{155537}{2^2 \cdot 3^2 \cdot  23 \cdot 101 \cdot 617}\eta_{11} +\frac{14306179}{3 \cdot 23 \cdot 101 \cdot 617}\eta_{12} + \frac{7867355933}{2^2 \cdot 3^2 \cdot  23 \cdot 101 \cdot 617}\eta_{13}, \\
    \omega_{11} &\coloneqq - \frac{71464984394573705}{2 \cdot 3^2 \cdot 23 \cdot 101 \cdot 617}\omega_5 - \frac{1741238355015323631}{2^3 \cdot 23 \cdot 101 \cdot 617}\omega_6 -5\eta_{11} -215\eta_{12} + 770\eta_{13}, 
    \end{align*}

\begin{align*}
    \omega_{12} &\coloneqq  \frac{5353742302066360397}{2^3 \cdot 3^2 \cdot 23 \cdot 101 \cdot 617}\omega_6 +\frac{952}{3^2 \cdot  23 \cdot 101 \cdot 617}\eta_8 
- \frac{90297452}{3^2 \cdot  23 \cdot 101 \cdot 617}\eta_{9} \\ &\hspace{2cm}+ \frac{90339725}{3^2 \cdot  23 \cdot 101 \cdot 617}\eta_{10} - \frac{7765582349}{2 \cdot 3^2 \cdot 23 \cdot 101 \cdot 617}\eta_{11} +
\frac{4635257585}{3 \cdot 23 \cdot 101 \cdot 617}\eta_{12} \\ &\hspace{2cm} + \frac{17104143219737}{2 \cdot 3^2 \cdot 23 \cdot 101 \cdot 617}\eta_{13}, \\
    \omega_{13} &\coloneqq -9\eta_7 + \frac{12912015}{23 \cdot 101 \cdot 617}\eta_8 - \frac{1109370333}{2 \cdot 23 \cdot 101 \cdot 617}\eta_9 +
\frac{2270264487}{2^2 \cdot  23 \cdot 101 \cdot 617}\eta_{10} \\ &\hspace{2cm} + \frac{15892479453}{2^3 \cdot 23 \cdot 101 \cdot 617}\eta_{11} + \frac{2443454506197}{2 \cdot 23 \cdot 101 \cdot 617}\eta_{12} +
\frac{397336735972815}{2^3 \cdot 23 \cdot 101 \cdot 617} \eta_{13}.
\end{align*}

Following $\S$\ref{subsec: calculating hodge}, for the matrix $Z = (\text{Tr}(T_3) I_{2g} - 2g T_3)C^{-1}$ we obtain

\begin{center}
$
    \begin{pmatrix}
    * & * & * & * & * & * & * &  14 &  0 &  14  & 0 &  0 &  0  & 0 \\
* & * & * & * & * & * & * &  -84 &  -56 &  -42 &  -28 &  -14 &  0  & 0 \\
* & * & * & * & * & * & *  & 42 &  0 &  14 &   0 &  0 &  0  & 14 \\
* & * & * & * & * & * & * &  364  & 168 &  168 &  70 &  84 &  42  & 14 \\
* & * & * & * & * & * & * &  0  & 56  & -14 &  28  & 14 &  0 &  -14 \\
    * & * & * & * & * & * & * &  -406 &   -126  & -196 &  -42 &  -84 &  -56 &  -28 \\
* & * & * & * & * & * & *  & 378 &  182 &  182 &  70 &  98  & 56 &  0 \\
-14  & 84 &  -42 &  -364  & 0  & 406 &  -378 &  0 &  0  & 0  & 0 &  0  & 0 &  0 \\
0  & 56 &  0  & -168  & -56 &  126 &  -182 &  0  & 0 &  0 &  0  & 0 &  0  & 0 \\
-14 &  42 &  -14 &  -168 &  14  & 196 &  -182 &  0 &  0  & 0  & 0 &  0  & 0  & 0 \\
0  & 28 &  0 &  -70  & -28 &  42 &  -70 &  0 &  0 &  0  & 0 &  0 &  0 &  0 \\
0  & 14 &  0 &  -84  & -14 &  84 &  -98 &  0 &  0 &  0 &  0 &  0 &  0 &  0 \\
0  & 0  & 0  & -42 &  0 &  56  & -56 &  0 &  0 &  0 &  0 &  0 &  0 &  0 \\
0  & 0  & -14  & -14 &  14  & 28 &  0 &  0 &  0 &  0 &  0 &  0 &  0 &  0 \\ 
\end{pmatrix},
$
\end{center}
where, for completeness, the top left $7 \times 7$ block is equal to 

\noindent\resizebox{\columnwidth}{!}{    
$
    \begin{pmatrix}
        0 & \frac{9540337662156}{23 \cdot 101 \cdot 617} & \frac{4335054124247763}{2 \cdot 23 \cdot 101 \cdot 617} & \frac{6401000762099027}{2 \cdot 3^2 \cdot 23 \cdot 101 \cdot 617} &  \frac{-13110192931669897}{2 \cdot 3 \cdot  23 \cdot 101 \cdot 617} &  \frac{-267727903087391908}{3^2 \cdot  23 \cdot 101 \cdot 617} &   \frac{-3048788107254906597}{2 \cdot 23 \cdot 101 \cdot 617} \\[6pt]
        \frac{-9540337662156}{23 \cdot 101 \cdot 617} &  0  & \frac{-11601405007652425}{2 \cdot 3^2 \cdot 23 \cdot 101 \cdot 617} &    \frac{-12300403690199765}{2 \cdot 3^2 \cdot 23 \cdot 101 \cdot 617}  & \frac{3796374812233423}{2 \cdot 3 \cdot  23 \cdot 101 \cdot 617} &  \frac{248689820059036475}{3^2 \cdot  23 \cdot 101 \cdot 617} &  \frac{6397851759663749133}{2^2 \cdot  23 \cdot 101 \cdot 617} \\[6pt]
        \frac{-4335054124247763}{2 \cdot 23 \cdot 101 \cdot 617} &  \frac{11601405007652425}{2 \cdot 3^2 \cdot 23 \cdot 101 \cdot 617} &  0 &   \frac{56474481585050101769}{2^2 \cdot 3^2 \cdot  23 \cdot 101 \cdot 617} &  \frac{54836979520019902345}{2^2 \cdot 3^2 \cdot  23 \cdot 101 \cdot 617} &    \frac{-279165429983233369247}{2^2 \cdot 3^2 \cdot  23 \cdot 101 \cdot 617} &  \frac{-6697486286619740202619}{2^2 \cdot 3^2 \cdot  23 \cdot 101 \cdot 617} \\[6pt]
        \frac{-6401000762099027}{2 \cdot 3^2 \cdot 23 \cdot 101 \cdot 617}  & \frac{12300403690199765}{2 \cdot 3^2 \cdot 23 \cdot 101 \cdot 617} &   \frac{-56474481585050101769}{2^2 \cdot 3^2 \cdot  23 \cdot 101 \cdot 617} &  0 &  \frac{19062442769792099321}{2 \cdot 3 \cdot  23 \cdot 101 \cdot 617} &   \frac{-167738977762662963811}{2^2 \cdot 3^2 \cdot  23 \cdot 101 \cdot 617}  & \frac{-2382715749575681741459}{2^2 \cdot 3 \cdot 23 \cdot 101 \cdot 617} \\[6pt]
        \frac{13110192931669897}{2 \cdot 3 \cdot  23 \cdot 101 \cdot 617} &  \frac{-3796374812233423}{2 \cdot 3 \cdot  23 \cdot 101 \cdot 617} &   \frac{-54836979520019902345}{2^2 \cdot 3^2 \cdot  23 \cdot 101 \cdot 617} &  \frac{-19062442769792099321}{2 \cdot 3 \cdot  23 \cdot 101 \cdot 617}  & 0 &   \frac{394979016632843924869}{2^2 \cdot 3^2 \cdot  23 \cdot 101 \cdot 617}  & \frac{1700048679479095053199}{3^2 \cdot  23 \cdot 101 \cdot 617} \\[6pt]
         \frac{267727903087391908}{3^2 \cdot  23 \cdot 101 \cdot 617}  & \frac{-248689820059036475}{3^2 \cdot  23 \cdot 101 \cdot 617} & \frac{279165429983233369247}{2^2 \cdot 3^2 \cdot  23 \cdot 101 \cdot 617}  & \frac{167738977762662963811}{2^2 \cdot 3^2 \cdot  23 \cdot 101 \cdot 617} &    \frac{-394979016632843924869}{2^2 \cdot 3^2 \cdot  23 \cdot 101 \cdot 617}  & 0  & \frac{14066222820396411352429}{2^2 \cdot 3^2 \cdot  23 \cdot 101 \cdot 617} \\[6pt]
         \frac{3048788107254906597}{2 \cdot 23 \cdot 101 \cdot 617} &  \frac{-6397851759663749133}{2^2 \cdot  23 \cdot 101 \cdot 617} &   \frac{6697486286619740202619}{2^2 \cdot 3^2 \cdot  23 \cdot 101 \cdot 617}  & \frac{2382715749575681741459}{2^2 \cdot 3 \cdot 23 \cdot 101 \cdot 617} &   \frac{-1700048679479095053199}{3^2 \cdot  23 \cdot 101 \cdot 617}  & \frac{-14066222820396411352429}{2^2 \cdot 3^2 \cdot  23 \cdot 101 \cdot 617} &  0
    \end{pmatrix}.
$
}

\par Next, we compute $\boldsymbol{\beta_{\text{Fil}}}, \gamma_{\text{Fil}}$ to complete the Hodge filtration computation. Recall that we are computing $b_{\text{Fil}} \in \Q_p^g$, and $\gamma_{\text{Fil
}} \in \text{span}\{\frac{s_{88}}{s_{90}}, \frac{s_{89}}{s_{90}}\}$ such that the following expression has no poles
\begin{equation*}
    c + \gamma_{\text{Fil}} - b_{\text{Fil}} A_2 - B_2^ \intercal A_2.
\end{equation*} 
For the expression $ c - B_2^ \intercal A_2$, we obtain, as a  $q$-expansion 
\begin{multline*}
     c - B_2^ \intercal A_2 = \frac{180506942}{3^3  \cdot 23  \cdot 101  \cdot 617}q^{-6} + \frac{72124682}{12899619}q^{-5} - \frac{175739501}{25799238}q^{-4} - \frac{1378473929}{2 \cdot 3^3  \cdot 23  \cdot 101  \cdot 617}q^{-3} \\ \hspace{2cm} - \frac{521048143}{12899619}q^{-2}  -
    \frac{1403673712}{12899619}q^{-1} + \frac{-5712282}{1433291} + O(q).
\end{multline*}
Writing the entries of the vector $A_2$ as rows of coefficients of the terms $q^{-9}$ to $q^{-1}$, we obtain
\begin{equation*}
    \begin{pmatrix}
        0 & 0 & 0 & 0 & 0 & 0 & 0 & 0 & -1 \\[6pt]
       0 & \frac{-34}{12899619}  & \frac{-34}{12899619} & \frac{-109}{12899619} &  \frac{-143}{12899619} &  \frac{-286}{12899619} &   \frac{-395}{12899619} &  
    \frac{-12900266}{12899619} &   \frac{-892}{12899619} \\[6pt]
    0 & 0 & 0 & 0 & 0 & 0 & -1 & 0 & 0 \\[6pt]
    0 & \frac{-3065}{12899619} & \frac{-3065}{12899619} &
    \frac{719497}{51598476} & \frac{707237}{51598476} & - \frac{25092001}{25799238} & \frac{2133971}{51598476} &
    \frac{3560705}{51598476} & \frac{2115581}{25799238} \\[6pt]
    0 & 0 & 0 & 0 & -1 & 0 & 0 & 0 & 0  \\[6pt]
    0 & \frac{119}{12899619} & \frac{-12899500}{12899619} & \frac{763}{25799238} & \frac{1001}{25799238} & \frac{1001}{12899619} & \frac{2765}{25799238} & \frac{4529}{25799238} & \frac{3122}{12899619} \\[6pt]
    -1 & \frac{3099}{2866582} & \frac{3099}{2866582} &
    \frac{-719061}{11466328} &  \frac{-706665}{11466328} &  \frac{-706665}{5733164} &  \frac{-2132391}{11466328} & 
    \frac{-3558117}{11466328} &  \frac{-2113797}{5733164}
    
    \end{pmatrix}
\end{equation*}
so that the top row corresponds to the first entry of $A_2$ which is $-q^{-1} + O(q)$. Writing the $q$-expansions of $\{\frac{s_{88}}{s_{90}}, \frac{s_{89}}{s_{90}}\}$ in the same way, we have 
\begin{equation*}
    \begin{pmatrix}
           0  &  1  &    1  &  -42  &  -41  &  -82  & -124 &  -207  & -245 \\
   0   &   0  &    0   &   1   &   1  &    2   &   3   &    5    &  6
    \end{pmatrix}.
\end{equation*}
Therefore, we obtain the following for the vector $b_{\text{Fil}}$
\begin{equation*}
\begin{pmatrix}
    \frac{1764687596}{3^2 \cdot  23 \cdot 101 \cdot 617} &  \frac{2465679139}{3^3  \cdot 23  \cdot 101  \cdot 617} &   \frac{2461515581}{2 \cdot 3^3  \cdot 23  \cdot 101  \cdot 617} &  
    \frac{1249246271}{2 \cdot 3^3  \cdot 23  \cdot 101  \cdot 617} &   \frac{-35867104}{3^3  \cdot 23  \cdot 101  \cdot 617} &   0 &   0
\end{pmatrix},
\end{equation*}
so that $\boldsymbol{\beta}_{\text{Fil}} = (0, ..., 0, b_\text{Fil})^\intercal$. For $\gamma_{\text{Fil}}$ we obtain
\begin{equation*}
    \gamma_{\text{Fil}} = \frac{148022444521}{2\cdot3^2 \cdot 23^2 \cdot 101^2 \cdot 617^2}\frac{s_{88}}{s_{90}} + \frac{75730523315027}{2^3 \cdot 23^2 \cdot 101^2 \cdot 617^2}\frac{s_{89}}{s_{90}},
\end{equation*}
and so writing $\gamma_{\text{Fil}}$ as a $q$-expansion we have 
\begin{multline*}
    \gamma_{\text{Fil}} = - \frac{148022444521}{36977815632258}q^{-8} - \frac{148022444521}{36977815632258}q^{-7} -
    \frac{218902313051905}{49303754176344}q^{-6} \\ \hspace{2cm} - \frac{657299028933799}{147911262529032}q^{-5} -
    \frac{657299028933799}{73955631264516}q^{-4} - \frac{1971304997023313}{147911262529032}q^{-3} \\ \hspace{0cm} 
    - \frac{365034551679203}{16434584725448}q^{-2} - \frac{1972193131690439}{73955631264516}q^{-1} + \\ \hspace{0cm}
    - \frac{46569944151137785}{73955631264516} - \frac{1972489176579481}{36977815632258}q -
    \frac{5916875439960359}{73955631264516}q^2 \\ - \frac{5041266187752739}{49303754176344}q^3 -
    \frac{20387519602731365}{147911262529032}q^4 + O(q^5).
\end{multline*}

Hence we have calculated the Hodge filtration ending with $\boldsymbol{\beta}_{\text{Fil}}^\intercal$ equal to

\noindent\resizebox{\columnwidth}{!}{    
$
    \begin{pmatrix}
        0 &   0& 0 &   0& 0 &   0& 0 &
    \frac{1764687596}{3^2 \cdot  23 \cdot 101 \cdot 617} &  \frac{2465679139}{3^3  \cdot 23  \cdot 101  \cdot 617} &  \frac{2461515581}{2 \cdot 3^3  \cdot 23  \cdot 101  \cdot 617} &  
    \frac{1249246271}{2 \cdot 3^3  \cdot 23  \cdot 101  \cdot 617} &   \frac{-35867104}{3^3  \cdot 23  \cdot 101  \cdot 617} &   0 &   0
    \end{pmatrix}.
$
}
and $\gamma_{\text{Fil}}$ equal to 
\begin{equation*}
\frac{148022444521}{2\cdot3^2 \cdot 23^2 \cdot 101^2 \cdot 617^2}\frac{s_{88}}{s_{90}} + \frac{75730523315027}{2^3 \cdot 23^2 \cdot 101^2 \cdot 617^2}\frac{s_{89}}{s_{90}}.
\end{equation*}

\bibliographystyle{alpha}
\bibliography{bibliography.bib}
\end{document}